\documentclass[11pt]{article}

\usepackage{amssymb,latexsym,amscd,array}
\usepackage{exscale}
\usepackage[centertags]{amsmath}
\usepackage{amssymb}
\usepackage{amsthm}
\usepackage{dsfont}
\usepackage[all]{xy}
\usepackage{lscape}
\usepackage{pdflscape}
\usepackage{color}
\usepackage{ upgreek }

\newtheorem{thm}{Theorem}[section]
\newtheorem{prop}[thm]{Proposition}

\newtheorem{rem}[thm]{Remark}
\newtheorem{lem}[thm]{Lemma}
\newtheorem{defn}[thm]{Definition}
\newtheorem{cor}[thm]{Corollary}

\newcommand{\mcd}{\mathcal{D}}

\newcommand{\mch}{\mathcal{H}}

\newcommand{\mco}{\mathcal{O}}

\newcommand{\mcr}{\mathcal{R}}

\newcommand{\mbc}{\mathbb{C}}
\newcommand{\mbd}{\mathbb{D}}

\newcommand{\mbg}{\mathbb{G}}

\newcommand{\mbn}{\mathbb{N}}

\newcommand{\mbp}{\mathbb{P}}

\newcommand{\mbr}{\mathbb{R}}

\newcommand{\mbv}{\mathbb{V}}

\newcommand{\p}{\partial}
\newcommand{\ra}{\rightarrow}
\newcommand{\lra}{\longrightarrow}
\newcommand{\FL}{\textup{FL}}

\numberwithin{equation}{subsection}

\begin{document}
\title{A comparison theorem between Radon and Fourier-Laplace transforms for D-modules}
\author{Thomas Reichelt}
\date{February 04, 2015}
\maketitle
\begin{abstract}
\noindent We prove a comparison theorem between the $d$-plane Radon transform and the Fourier-Laplace transform for $D$-modules. This generalizes results of Brylinski and d'Agnolo-Eastwood.
\end{abstract}
\renewcommand{\thefootnote}{}
\footnote{
During the preparation of this paper, the author was supported by a postdoctoral fellowship of the ``Fondation sciences math\'ematiques de Paris'',
by the DFG grant He 2287/2-2
and received partial support by the ANR grant ANR-08-BLAN-0317-01 (SEDIGA).
}

\section*{Introduction}
The history of the Radon transform goes back to the famous paper \cite{Radon} of Radon. This transformation associates to a function $f$ on $\mbr^n$ a corresponding function on the family of affine $d$-dimensional planes in $\mbr^n$ whose value at a given plane is the integral over $f$ restricted to this plane.  Since then various generalizations, like the Radon transform on homogeneous spaces \cite{Helga1} and the Penrose transform \cite{EPW} were made, which had plenty of applications in representation theory, harmonic analysis and mathematical physics. Later it was realized that the Radon transform and its various cousins could be best understood in the context of integral geometry and $\mcd$-modules (see \cite{Agno1} for a nice overview).  This Radon transform for $\mcd$-modules was introduced by Brylinski in \cite{Brylinski}. There, he considers $\mcd$-modules on the complex projective space and measures their restriction to all $d$-planes, which gives rise to a (complex of) $\mcd$-module(s) on the corresponding Grassmanian. As an application he proved, among other things, the irreducibility of the monodromy action on the vanishing cohomology of a hyperplane section of a (possibly singular) variety. The study of these integral transforms was further advanced by d'Agnolo and Schapira who on the one hand extended in \cite{dAS1} the hyperplane  Radon transform to twisted $\mcd$-modules and to a quantized contact transformation between the cotangent bundles of the projective spaces and on the other hand to the very general setting of double fibrations in \cite{dAS2}.\\ 

The Fourier-Laplace transform on the other hand is an indispensable tool in the theory of differential equations. In the context of integral geometry it is a transform with a so-called exponential kernel. This property is reflected in the $\mcd$-module picture by the fact that the Fourier-Laplace transform does not preserve regular holonomicity. If however the $\mcd$-module is monodromic then the Fourier-Laplace transform is equivalent to the so-called Fourier-Sato transform (or monodromic Fourier-Laplace transform) which preserves regularity. By a theorem of Brylinski \cite{Brylinski},  the hyperplane Radon transform for $\mcd$-modules on $\mbp^n$ is closely related to the monodromic Fourier-Laplace transform on $\mbc^{n+1}$. Roughly speaking the theorem of Brylinski says, that the hyperplane Radon transform of a holonomic $\mcd$-module on $\mbp^n$ is isomorphic to the (monodromic) Fourier-Laplace transform of a specific lift of the $\mcd$-module from $\mbp^n$ to $\mbc^{n+1}$. This result was generalized by d'Agnolo and Eastwood \cite{AE} to quasi-coherent $\mcd$-modules and to a variant of the Radon transform which does measure the restriction of the $\mcd$-module to the complement of a given hyperplane, rather then the restriciton to the hyperplane itself. \\

In this paper we extend the results of Brylinski and d'Agnolo-Eastwood to the case of the $d$-plane Radon transform. Possible applications of this result are  explicit computations of Gauss-Manin systems of maps whose fibers are $d$-plane sections inside a quasi-projective variety. This technique was already used in \cite{Reich2} in the case of families of Laurent polynomials, where the fiber of a given Laurent polynomial can be compactified to a hyperplane section of a projective toric variety, which is given by the Newton polytope of the corresponding Laurent polynomial. In forthcoming work it is planned to use the comparison theorem between the $d$-plane Radon transform and the Fourier-Laplace transform, which is proven here, to study hyperplane arrangements. More precisely, we want to use a variant of the $d$-plane Radon transform, defined below, to compute explicitly the Gauss-Manin system of the universal family of hyperplane arrangements and of the universal family of the complements. \\

Let us give a short overview of the paper. We first give a brief review of algebraic $\mcd$-modules, the Radon as well as the Fourier-Laplace transform and state the comparison result for the hyperplane case of d'Agnolo and Eastwood. The proof of the general case proceeds as follows. We first prove that the $d$-plane Radon transform for a $\mcd$-module is equivalent to a diagonal embedding of this $\mcd$-module on a product of projective spaces and then applying the hyperplane Radon transform on each factor (Lemma \ref{lem:RadontoRadonprod} and Lemma \ref{lem:RadonprodtoRadonpart}). This result stems from the simple geometric fact that a $d$-plane in $\mbp^n$ is isomorphic to the diagonal in $(\mbp^n)^{\times n-d}$ cut with an appropiate hyperplane on each factor. We then prove in Proposition \ref{prop:RadonparttoFLpart} that we can apply the theorem of d'Agnolo and Eastwood on each factor. Finally we have to prove that the extension functors and the Fourier-Laplace transform interchange (Lemma  \ref{lem:FlparttoFLprod}).

\textbf{Acknowledgements:}  I would like to thank Christian Sevenheck and Claude Sabbah for useful discussions. Furthermore, I thank Claus Hertling for his continuous support and interest in my work.

\section{Preliminaries}

In the first section we review briefly the theory of algebraic $\mcd$-modules in order to fix notations. In the second section we review the $d$-plane Radon transform and various variants (cf. Definition \ref{def:radtrafo}) and show how they arise as an integral transformation (cf. Proposition \ref{prop:kernelRadcompare}). Then we show in Proposition \ref{prop:Radpreserve} that the various transformations preserve (regular) holonomicity and \mbox{(quasi-)coherence}. In the last section we introduce in Definition \ref{def:FLtrafo} the Fourier-Laplace transform for $\mcd$-modules as an integral transform with exponential kernel .

\subsection{$\mcd$-modules}
Let $X$ be a smooth algebraic variety over $\mbc$. We denote by $M(\mcd_X)$ the abelian category of algebraic left $\mcd_X$-modules on $X$. The full triangulated subcategories of $D^b(\mcd_X) := D^b(M(\mcd_X))$, consisting of objects with $\mco_X$-quasi-coherent resp. $\mcd_X$-coherent resp. (regular) holonomic cohomology are denoted by $D^b_{qc}(\mcd_X)$ resp. $D^b_{coh}(\mcd_X)$ resp. $D^b_{h}(\mcd_X)$ resp. $D^b_{rh}(\mcd_X)$. We subsume the different cases by writing $D^b_\ast(\mcd_X)$ for $\ast \in \{qc, coh, h ,rh\}$.

Let $f: X \ra Y$ be a map between smooth algebraic varieties. We denote by $\mcd_{X \ra Y}$ resp. $\mcd_{Y \leftarrow X}$ the transfer bimodules. Let $M \in D^b(\mcd_X)$ and $N \in D^b(\mcd_Y)$, then the direct and inverse image for $\mcd$-modules is defined by
\begin{align}
f_+M &:= Rf_* (\mcd_{Y \leftarrow X} \overset{L}{\otimes}_{\mcd_X} M) , \notag \\
f^+ N &:= \mcd_{X \ra Y} \overset{L}{\otimes}_{f^{-1} \mcd_Y} f^{-1}N\, . \notag
\end{align}
Recall that the functors $f_+,f^+$ preserve quasi-coherence, holonomicity and regular holonomicity (see e.g., \cite{Hotta}).

If $f:X \ra Y$ is non-characteristic, then the functor $f^+$ preserves coherency and is exact.\\

\noindent Denote by $\omega_X$ the canonical line bundle on $X$. There exists a duality functor $\mbd: D^b_{coh}(\mcd_X) \ra D^b_{coh}(\mcd_X)$ defined by
\[
\mbd M := \mch om_{\mco_X}(\omega_X,R\mch om_{\mcd_X}(M,\mcd_X))[dim X].
\]
Recall that for a single holonomic $\mcd_X$-module $M$, the holonomic dual is also a single holonomic $\mcd_X$-module (\cite[Corollary 2.6.8 (iii)]{Hotta}).

For a morphism $f: X \ra Y$ between smooth algebraic varieties we additionally define the functors $f_\dag := \mbd \circ f_+ \circ \mbd$ and $f^\dag := \mbd \circ f^+ \circ \mbd$.

Consider the following cartesian diagram of algebraic varieties
\[
\begin{xy}
\xymatrix{ Z \ar[r]^{f'} \ar[d]_{g'}  & W \ar[d]^g \\ Y \ar[r]^f &  X}
\end{xy}
\]
then we have the canonical isomorphism $f^+ g_+[d] \simeq  g'_+ f{'}^+[d']$, where $d:= \dim Y - \dim X$ and $d' := \dim Z - \dim W$(cf. \cite[Theorem 1.7.3]{Hotta}).\\
 
Notice that by symmetry we have also the canonical isomorphism  $g^+ f_+[\tilde{d}] \simeq f'_+ g{'}^+[\tilde{d}']$ with $\tilde{d}:= \dim W - \dim X$ and $\tilde{d}':= \dim Z - \dim Y$. In the former
case we say we are doing a base change with respect to $f$, in the latter case with respect to $g$.\\

´
Using the duality functor, we get isomorphisms:
\[
f^\dag g_\dag[-d] \simeq g'_\dag f{'}^\dag[-d'] \qquad \text{and} \qquad g^\dag f_\dag[-\tilde{d}] \simeq f'_{\dag}\,g{'}^\dag[-\tilde{d}']\, .
\]

Let $M \in D^b(\mcd_X)$ and $N \in D^b(\mcd_Y)$. We denote by $M\boxtimes N \in D^b(\mcd_{X \times Y})$ the exterior tensor product. The exterior tensor product preserves quasi-coherence, coherence, holonomicity and regular holonomicity. If  $M_1,M_2 \in D^b(\mcd_X)$, we denote by
\[
M_1 \overset{L}{\otimes} M_2 := \Delta^+ (M_1 \boxtimes M_2)
\]
the internal tensor product, where $\Delta: X \ra X\times X$ is the diagonal embedding. The internal tensor product preserves quasi-coherence, holonomicity and regular holonomicity. Notice that it preserves coherency if $M_1 \boxtimes M_2$ is noncharacteristic with respect to $\Delta$, i.e. if $M_1, M_2$ satisfy the following transversality condition
\[
char(M_1) \cap  char(M_2) \subset T^*_{X} X\, .
\]

Let $f:X \ra Y$ a map between smooth algebraic varieties. One has the following  projection formula (cf. \cite[Corollary 1.7.5]{Hotta})
\begin{equation}\label{eq:projformula}
f_+(M \overset{L}{\otimes} f^+ N) \simeq f_+ M \overset{L}{\otimes} N\, .
\end{equation}

In the following it will sometimes be convenient to use the language of integral kernels. Let $X$ and $Y$ be two smooth varieties, $M \in D^b(\mcd_X)$ and $K \in D^b(\mcd_{X \times Y})$. 
Denote by $q_1: X \times Y \ra X$ resp. $q_2: X \times Y \ra Y$ the projection to the first resp. second factor.
The integral transform with respect to the kernel $K$ is defined by
\begin{align}
\diamond\, K : D^b(\mcd_X) &\lra D^b(\mcd_Y)\, , \notag \\
M &\mapsto M \diamond K = q_{2+}(q_1^+M \overset{L}{\otimes} K)\, .\notag
\end{align}

Let $Z$ be another smooth algebraic variety and $\widetilde{K} \in D^b(\mcd_{Y \times Z})$. The convolution of the two kernels $K$ and $\widetilde{K}$ is defined by
\[
K \diamond \widetilde{K} := q_{13 +}(q_{12}^+ K \overset{L}{\otimes} q_{23}^+ \widetilde{K})\, ,
\]
where $q_{ij}$ is the projection from $X \times Y \times Z$ to the corresponding factor. Notice that the folding is associative in the sense that for $M \in D^b(\mcd_X)$ we have
\[
(M \diamond K) \diamond \widetilde{K} \simeq M \diamond ( K \diamond \widetilde{K}) \in D^b(\mcd_Z)\, .
\]

\subsection{Radon transform}
Let $V$ be a complex $n+1$-dimensional vector space and denote by $\hat{V}$ the dual vector space. Fix $k \in \{1, \ldots , n\}$ and set $\hat{W}:=\hat{V}^{k}$ , so a point in $\hat{W}$ gives rise to a vector subspace of $V$ of codimension $\leq k$. Set $d:= n-k$ and denote by $S(k,n)$ the subvariety of $\hat{W}$ consisting of points giving rise to vector subspaces of codimension equal to $k$ resp. $d$-planes in $\mbp(V)$. Notice that $GL(k)$ acts on $S(k,n)$ from the left. The quotient
\[
G(d,n):= Gl(k) \setminus S(k,n)
\]
is the Grassmanian parametrizing $n+1-k$-dimensional subspaces in $V$, i.e. $d$-planes in $\mbp(V)$. We use the following abbreviations $\mbp:= \mbp(V)$, $\dot{V} := V \setminus \{ 0\}$ and $\mbg := G(d,n)$.

Let $Z \overset{i_Z}{\hookrightarrow} \mbp \times \mbg$ be the universal hyperplane, i.e.
\[
Z:= \{ [v], [\lambda^1, \ldots, \lambda^k] \in \mbp \times \mbg \mid \lambda^1(v) = \ldots = \lambda^k(v) = 0\}
\]
and denote by $C \overset{j_C}{\hookrightarrow} \mbp \times \mbg$ its complement

\[
C:= \{ [v], [\lambda^1, \ldots, \lambda^k] \in \mbp \times \mbg \mid \exists l \in \{1, \ldots ,k\}\;  \text{with}\; \lambda^l(v) \neq 0\}\, .
\]

We will define various versions of the Radon transform. Consider the following diagram
$$
\xymatrix{ && C \ar[drr]^{\pi_2^C} \ar[dll]_{\pi_1^C} \ar@{^(->}[d]^{j_C}&& \\ \mbp && \mbp \times \mbg \ar[ll]_{\pi_1} \ar[rr]^{\pi_2} && \mbg \; , \\ && Z \ar[ull]^{\pi_1^Z} \ar@{^(->}[u]_{i_Z} \ar[rru]_{\pi_2^Z} &&}
$$
where $\pi_1$ resp. $\pi_2$ are the projections to the first resp. second factor and $\pi_1^C, \pi_1^Z, \pi_2^C, \pi_2^Z$ are the corresponding restrictions of $\pi_1$ and $\pi_2$ to $C$ resp. $Z$.

\begin{defn}\label{def:radtrafo}
Let $V, Z , U$ as above. The Radon transform is the functor 
\begin{align}
\mcr_\delta: D^b_{qc}(\mcd_{\mbp}) &\ra D^b(\mcd_{\mbg}) \notag \\
M &\mapsto \pi^Z_{2+}(\pi^Z_1)^+ M \simeq \pi_{2+}i_{Z+}i^+_Z\pi^+_1 M\, .\notag
\end{align}
Define two variants
\begin{align}
\mcr_{1/t}: D^b_{qc}(\mcd_{\mbp}) &\ra D^b(\mcd_{\mbg}) \notag \\
M &\mapsto \pi^C_{2+}(\pi^C_1)^+ M \simeq \pi_{2+}j_{C+}j^+_C\pi^+_1 M \notag
\end{align}
and
\begin{align}
\mcr_Y: D^b_{coh}(\mcd_{\mbp}) &\ra D^b(\mcd_{\mbg}) \notag \\
M &\mapsto \pi^C_{2\dag}(\pi^C_1)^+ M \simeq \pi_{2+}j_{C\dag}j^+_C\pi^+_1 M\, , \notag
\end{align}
as well as a constant Radon transform 
\begin{align}
\mcr_1: D^b_{qc}(\mcd_{\mbp}) &\ra D^b(\mcd_{\mbg}) \notag \\ M &\mapsto \pi_{2+}(\pi_1)^+ M\, . \notag
\end{align}
\end{defn}

In order to compare the Radon transform to the Fourier-Laplace transform which will be introduced below, we will need another type of Radon transform but this time with target in $\hat{W} = \hat{V}^{\times k}$ instead of $\mbg = G(n-k,n)$.

Let $Z' \overset{i_{Z'}}{\hookrightarrow} \mbp \times \hat{W}$ be defined as
\[
Z':= \{ [v], \lambda^1, \ldots, \lambda^k \in \mbp \times \hat{W} \mid \lambda^1(v) = \ldots = \lambda^k(v) = 0\}
\]

and denote by $C' \overset{j_{C'}}{\hookrightarrow} \mbp \times \hat{W}$ the complement of $Z'$:
\[
C':= \{ [v], \lambda^1, \ldots, \lambda^k \in \mbp \times \hat{W} \mid \exists l \in \{1, \ldots ,k \} \; \text{with}\; \lambda^l(v) \neq 0 \}\, .
\]
We will also consider the subvariety $A'\overset{i_{A'}}{\hookrightarrow}\mbp \times \hat{W}$:
\[
A' := \{ [v], \lambda^1, \ldots, \lambda^k \in \mbp \times \hat{W} \mid \exists l \in \{1, \ldots , k\} \; \text{with}\; \lambda^l(v) = 0\}\, 
\]
together with its complement $U' \overset{j_{U'}}{\lra} \mbp \times \hat{W}$:
\[
U' := \{ [v], \lambda^1, \ldots, \lambda^k \in \mbp \times \hat{W} \mid \lambda^1(v)\neq 0 , \ldots ,\lambda^k(v) \neq 0\}\, .
\]
Consider the following diagrams
$$
\xymatrix{ && C' \ar[drr]^{\pi_2^{C'}} \ar[dll]_{\pi_1^{C'}} \ar@{^(->}[d]^{j_{C'}}&& \\ \mbp && \mbp \times \hat{W} \ar[ll]_{\pi_1} \ar[rr]^{\pi_2} && \hat{W} \\ && Z' \ar[ull]^{\pi_1^{Z'}} \ar@{^(->}[u]_{i_{Z'}} \ar[rru]_{\pi_2^{Z'}} &&} \phantom{ttttttttt} \xymatrix{ && U' \ar[drr]^{\pi_2^{U'}} \ar[dll]_{\pi_1^{U'}} \ar@{^(->}[d]^{j_{U'}}&& \\ \mbp && \mbp \times \hat{W} \ar[ll]_{\pi_1} \ar[rr]^{\pi_2} && \hat{W}\; , \\ && A' \ar[ull]^{\pi_1^{A'}} \ar@{^(->}[u]_{i_{A'}} \ar[rru]_{\pi_2^{A'}} &&}
$$
where as above $\pi_1$ resp. $\pi_2$ are the projections to the first resp. second factor and $\pi_1^{C'}, \pi_1^{Z'},\pi_1^{U'},\pi_1^{A'}$, $\pi_2^{C'}, \pi_2^{Z'},\pi_2^{U'},\pi_2^{A'}$ are the corresponding restrictions of $\pi_1$ and $\pi_2$ to $C'$ resp. $Z'$ resp. $U'$ resp. $A$.
\begin{defn}
Let $Z' , U', C'$ as above. The affine Radon transform is the functor 
\begin{align}
\mcr'_\delta: D^b_{qc}(\mcd_{\mbp}) &\ra D^b(\mcd_{\hat{W}}) \notag \\
M &\mapsto \pi^{Z'}_{2+}(\pi^{Z'}_1)^+ M \simeq \pi_{2+}i_{Z'+}i^+_{Z'} \pi^+_1 M\, . \notag 
\end{align}
Define the variants
\begin{align}
\mcr'_{1/t}: D^b_{qc}(\mcd_{\mbp}) &\ra D^b(\mcd_{\hat{W}}) \notag \\
M &\mapsto \pi^{C'}_{2+}(\pi^{C'}_1)^+ M \simeq \pi_{2+}j_{C'+}j^+_{C'}\pi^+_1 M \notag\\
\mcr'_Y: D^b_{coh}(\mcd_{\mbp}) &\ra D^b(\mcd_{\hat{W}}) \notag \\
M &\mapsto \pi^{C'}_{2\dag}(\pi^{C'}_1)^+ M \simeq \pi_{2+}j_{C'\dag}j^+_{C'}\pi^+_1 M \notag \\
\mcr'_{U^+}: D^b_{qc}(\mcd_{\mbp}) &\ra D^b(\mcd_{\hat{W}}) \notag \\
M &\mapsto \pi^{U'}_{2+}(\pi^{U'}_1)^+ M \simeq \pi_{2+}j_{U'+}j^+_{U'}\pi^+_1 M \notag
%\mcr'_{U^\dag}: D^b_{coh}(\mcd_{\mbp}) &\ra D^b(\mcd_{\hat{W}}) \notag \\
%M &\mapsto \pi^{U'}_{2\dag}(\pi^{U'}_1)^+ M \simeq \pi_{2+}j_{U'\dag}j^+_{U'}\pi^+_1 M\, , \notag
\end{align}
as well as a constant Radon transform 
\begin{align}
\mcr'_1: D^b_{qc}(\mcd_{\mbp}) &\ra D^b(\mcd_{\hat{W}}) \notag \\
M &\mapsto \pi_{2+}(\pi_1)^+ M\, . \notag
\end{align}
\end{defn}

We can now express the affine Radon transforms $\mcr'_u$ for $u\in \{1,\delta, 1/t, Y, U^+\}$ as an integral transformation with respect to a kernel $R'_u$. Set
\[
R'_\delta = i_{Z' +} \mco_{Z'}, \;\, R'_Y = j_{C' \dag}\mco_{C'}, \;\, R'_{1/t} = j_{C' +} \mco_{C'}, \;\, R'_{U^+} = j_{U' +} \mco_{U'}, \;\,  R'_1 = \mco_{\mbp \times \hat{W}}\, .
\]
Recall that
\[
M \diamond R'_u = \pi_{2+}(\pi_1^+(M) \overset{L}{\otimes} R'_u)\, ,
\]
then we have the following result.
\begin{prop}\label{prop:kernelRadcompare} Let $M \in D^b_{qc}(\mbp)$. We have the following isomorphism
\[
M \diamond R'_{\hat{u}} \simeq \mcr'_{\hat{u}}(M) \quad \text{for} \; \hat{u} \in \{ \delta,1, 1/t, U^+ \}\, ,
\]
if $M \in D^b_{coh}(\mbp)$ then
\[
M \diamond R'_{\hat{u}} \simeq \mcr'_{\hat{u}}(M)\quad \text{for} \; \hat{u} \in \{Y \}\,.
\]
\end{prop}

\begin{proof}
The proof essentially uses the projection formula (cf. \eqref{eq:projformula}). We prove the statement for $u = \delta$ the other cases are similar or easier.
\begin{align}
M \diamond R'_\delta &= \pi_{2 +}(\pi_1^+(M) \overset{L}{\otimes} i_{Z' +} \mco_{Z'} ) \notag \\
&\simeq \pi_{2+} i_{Z' +} (i_{Z'}^+ \pi_1^+M \overset{L}{\otimes} \mco_{Z'}) \notag \\
&\simeq \pi_{2+}^{Z'}((\pi_1^{Z'})^+ M) \notag \\
&= \mcr'_\delta(M) \notag
\end{align}
\end{proof}

Finally, we define the last variant with respect to a kernel having support on the non-smooth subvariety $A'$.
\begin{defn}
Let $A' \subset \mbp \times \hat{W}$ as above. Set
\[
R'_{A^+} := R\Gamma_{A'} \mco_{\mbp \times \hat{W}}% \quad \text{and} \quad R'_{A^\dag} := \mbd R'_{A^+}
\]
and define for $M \in D^b_{qc}(\mbp)$
\[
\mcr'_{A^+} (M) := \pi_{2+}(\pi_1^+(M)\overset{L}{\otimes} R'_{A^+})\simeq \pi_{2+} R\Gamma_{A'}\, \pi^+_1M\, . %\quad \text{and} \quad \mcr'_{A^\dag} (M) := \pi_{2+}(\pi_1^+(M)\overset{L}{\otimes} R'_{A^\dag})\, .
\]
\end{defn}

\begin{prop}\label{prop:Radpreserve}
The Radon transforms preserve the following subcategories 
\begin{align}
\mcr_\delta, \mcr_{1/t}, \mcr_1 : D^b_*(\mcd_\mbp) &\ra D^b_*(\mcd_\mbg) \quad \text{for}\; * \in \{qc,coh,h,rh\}\, , \notag \\
\mcr_Y : D^b_*(\mcd_\mbp) &\ra D^b_*(\mcd_\mbg) \quad \text{for}\; * \in \{coh,h,rh\}\, , \notag
\end{align}
The affine Radon transforms preserve the following subcategories
\begin{align}
\mcr'_\delta, \mcr'_{1/t}, \mcr_{U^+}, \mcr'_1, \mcr'_{A^+} : D^b_*(\mcd_\mbp) &\ra D^b_*(\mcd_{\hat{W}}) \quad \text{for}\; * \in \{qc,coh,h,rh\}\, , \notag \\
\mcr'_Y : D^b_*(\mcd_\mbp) &\ra D^b_*(\mcd_{\hat{W}}) \quad \text{for}\; * \in \{coh,h,rh\}\, , \notag
\end{align}
\end{prop}
\begin{proof}
First notice that the claim is clear for $* \in \{qc,h,rh\}$, because the direct image, the  inverse image and the derived tensor product $\overset{L}{\otimes}$ preserve quasi-coherence, holonomicity and regular holonomicity and the proper direct image preserves holonomicity and regular holonomicity. The functors $\mcr_\delta, \mcr_{Y}, \mcr_1$ resp. $\mcr'_\delta, \mcr'_{Y}, \mcr'_1$ preserve coherency because $\pi^Z_1, \pi^C_1$ resp. $\pi^{Z'}_1, \pi^{C'}_1$ are smooth and coherency is preserved by proper direct images.  In order to prove that $\mcr'_{1/t}$ preserves coherency, recall that we have the isomorphism
\[
\mcr'_{1/t}(M) \simeq M \diamond R'_{1/t} \simeq \pi_{2+}(\pi_1^+M \overset{L}{\otimes} R'_{1/t})\, .
\]
Because $\pi_2$ is proper, it is enough to show that $\pi_1^+ M \overset{L}{\otimes} R'_{1/t} \in D^b_{coh}(\mcd_{\mbp \times \hat{W}})$ for $M \in D^b_{coh}(\mcd_{\mbp})$. Notice that $char(R'_{1/t}) \subset T^*_{Z'}\mbp \times \hat{W}\, \cup\, T^*_{\mbp \times \hat{W}} \mbp \times \hat{W}$. Thus one can easily compute that the transversality condition 
\[
(char(M) \times T^*_{\hat{W}} \hat{W}) \cap \left(T^*_{Z'} (\mbp \times \hat{W}) \cup T^*_{\mbp \times \hat{W}} \mbp \times \hat{W} \right) \subset T^*_{\mbp \times \hat{W}} \mbp \times \hat{W}
\]
is satisfied. This shows the claim for $\mcr'_{1/t}$.
 The proofs for $\mcr'_{U^+},\mcr'_{A^+}$ and $\mcr_{1/t}$ can be easily adapted.
\end{proof}

\begin{lem}
Let $M \in D^b_{coh}(\mcd_\mbp)$, we have the following triangles in $D^b_{coh}(\mcd_{\hat{W}})$:
\begin{enumerate}
\item $\mcr'_1(M) \lra \mcr'_{1/t}(M) \lra \mcr'_\delta(M) \overset{+1}{\lra}$
\item $\mcr'_\delta(M) \lra \mcr'_Y(M) \lra \mcr'_1(M) \overset{+1}{\lra }$
\item $ \mcr'_{A^+}(M) \lra \mcr'_1(M) \lra \mcr'_{U^+}(M) \overset{+1}{\lra}$
%\item $\mcr'_{U^\dag}(M) \lra \mcr'_1(M) \lra \mcr'_{A^\dag}(M)  \overset{+1}{\lra }$
\end{enumerate}
\end{lem}
\begin{proof}
The first and the third triangle can be deduced from the following triangles
\[
i_{Z' +} i_{Z'}^+[-1] (N) \lra  N \lra j_{C' + } j_{C'}^+(N) \overset{+1}{\lra} \quad\text{and}\quad R\Gamma_A (N) \lra N \lra j_{U'+} j_{U'}^+(N) \overset{+1}{\lra}
\]
where $N \in D^b_{coh}(\mbp \times \hat{W})$. The second triangle can be deduced by applying $\mbd$ to the first one.
\end{proof}

\subsection{Fourier-Laplace transform} 
In the next definition we want to introduce the Fourier-Laplace transform.
\begin{defn}\label{def:FLtrafo}
Let as above $V$ be a vector space and denote by $\hat{V}$ its dual. Define $L := \mco_{V \times \hat{V}} e^{-\langle v, \hat{v} \rangle}$ where $\langle \bullet, \bullet \rangle$ is the natural pairing between $V$ and $\hat{V}$ and the $\mcd$-module structure is given by the product rule. Denote by  $\pi_1: V \times \hat{V} \ra V$, $\pi_2 : V \times \hat{V} \ra \hat{V}$ the canonical projections. The Fourier-Laplace transform is then defined by
\begin{align}
FL: D^b_{qc}(\mcd_V) &\lra D^b_{qc}(\mcd_{\hat{V}})\, , \notag \\
M &\mapsto \pi_{2 +}(\pi_1^+ M \overset{L}{\otimes} L)\, . \notag
\end{align}
\end{defn}
\begin{rem}
In the setting above with $\hat{W}= \hat{V}^k$ and $W:= V^k$ we can perform the Fourier-Laplace transform stepwise. Set $N_{j,l}:= V^{\times j} \times \hat{V}^{\times l}$ for $j,l \in \{0, \ldots , k\}$ with $j+l =k+1$ and denote by $v^1, \ldots , v^j, \lambda^1, \ldots \lambda^{l}$ elements of $N_{j,l}$. Define $L_{jl} := \mco_{N_{j,l}} e^{-\langle v^j,\lambda^1 \rangle}$. Notice that $L_{jl} \simeq (\pi^{j,l})^+ L$ if $\pi^{j,l}: N_{j,l} \ra V \times \hat{V}$ is the canonical projection to the $j-th$ factor of $V^{\times j}$ and to the first factor of $\hat{V}^{\times l}$. Set
\[
\FL_j (M) := (\pi_2^{j,l})_+((\pi_1^{j,l})^+M \otimes L_{jl})\, ,
\]
where
\begin{align}
\pi_1^{j,l}: N_{j,l} &\lra N_{j,l-1}\, , \notag \\
(v^1, \ldots , v^j, \lambda_1, \ldots , \lambda^l)&\mapsto (v^1, \ldots ,v^j, \lambda^2, \ldots, \lambda^l)\, , \notag \\
\pi_2^{j,l}: N_{j,l} &\lra N_{j-1,l}\, , \notag \\
(v^1, \ldots , v^j, \lambda_1, \ldots , \lambda^l)&\mapsto (v^1, \ldots ,v^{j-1}, \lambda^1, \ldots, \lambda^l)\, . \notag
\end{align}
One can easily show that for $M \in D^b_{qc}(\mcd_W)$ we have
\[
\FL_1 \circ \ldots \circ \FL_k(M) \simeq \FL(M)\, ,
\]
where the Fourier-Laplace transform on the right hand side is performed with respect to $W$.
\end{rem}

We will also need a relative version of the Fourier-Laplace transform.
So let $E \ra X$ be a vector bundle, $\hat{E} \ra X$ its dual and denote by $can: E \times_X \hat{E} \ra X$ the canonical pairing with respect to the fibers. Define $L_X := \mco_{E \times  \hat{E}}\, e^{-can}$. The Fourier-Laplace transform with respect to the base $X$ is defined by
 \begin{align}
 FL_X: D^b_{qc}(\mcd_E) &\lra D^b_{qc}(\mcd_{\hat{E}}) \notag \\
 M &\mapsto \pi^X_{2+} (( \pi_1^X)^+M \overset{L}{\otimes} L_X )\notag
 \end{align}
 where $\pi_1^X : E \times_X \hat{E} \ra E$ and $\pi_2^X: E \times_X \hat{E} \ra \hat{E}$ are the canonical projections.
 
Consider the case where the vector bundle $E$ is trivial, i.e. $\pi: E = X \times V \ra X$ and $\hat{E} = X \times \hat{V}$. Let $v_1, \ldots , v_n$ be coordinates on $V$ and denote by $\lambda_1, \ldots , \lambda_n$ the dual coordinates on $\hat{V}$. This gives an isomorphism between $E$ and $\hat{E}$ which allows us to view $\pi^{-1}\mcd_{X}$-modules to be defined on $\hat{E}$. Let $M$ be a $\mcd_{X \times V}$-module. We denote by $\hat{M}$ the $D_{X \times \hat{V}}$-module on $\hat{E}$ which is equal to $M$ as a $\pi^{-1}\mcd_X$-module and where $\lambda_i$ acts as $\p_{v_i}$ and $\p_{\lambda_i}$ acts as $-v_i$.

The following proposition goes back to Laumon and Katz \cite{KaLau}. The $\mcd$-module case is due to Malgrange \cite{Mal9}.
\begin{prop}\label{prop:FLconc}
There is the following isomorphism in $D^b(\mcd_{\hat{E}})$:
\[
FL_X(M) \simeq \hat{M} \qquad \text{for} \quad M \in M_{qc}(\mcd_{E}) \, .
\]
\end{prop}
\begin{proof}
The proof given in \cite{Mal9} of the case $X = pt$ carries over almost word for word.
\end{proof}

\section{The comparison theorem}

In this section we give the proof of the comparison theorem. In the first section we review the theorem of d'Agnolo and Eastwood for the codimension one case. In the second section we state the results for the higher codimension case and in the last section we give the proof of these statements.
\subsection{The case of codimension one}

Notice that in codimension one the varieties  $Z'$ and $A'$ as well as $U'$ and $C'$ are equal, hence we have the following equality of functors $\mcr_\delta[-1] = \mcr_{A^+}$ and $\mcr_{1/t} = \mcr_{U^+}$.

Denote by $\widetilde{\pi}: \widetilde{V} \lra V$ the blowup of the origin $0$ in $V$ and by $E$ the exceptional divisor. Then $\widetilde{V}$ carries the following stratification
\[
E \overset{i_E}{\hookrightarrow} \widetilde{V} \overset{j_{\dot{V}}}{\hookleftarrow} \dot{V}\, .
\]
Denote by $\widetilde{i}: \widetilde{V} \ra \mbp \times V$ the natural embedding. Consider the maps
\begin{equation}\label{eq:blowupmaps}
\mbp \overset{\widetilde{\pi}}{\leftarrow} \widetilde{V} \overset{\widetilde{j}}{\rightarrow} V
\end{equation}
obtained by restriction of the natural projections from $\mbp \times V$ and by
\[
\mbp \overset{\pi}{\leftarrow} \dot{V} \overset{j}{\rightarrow} V
\]
their restriction to $\dot{V}$. This gives rise to the following diagram
\begin{equation}\label{eq:diagram}
\xymatrix{ && \dot{V} \ar[d]^{j_{\dot{V}}}\ar[drr]^j \ar[dll]_\pi && \\
\mbp && \widetilde{V} \ar[rr]^{\widetilde{j}} \ar[ll]_{\widetilde{\pi}} && V\, . \\
&&  E \ar[u]_{i_E} \ar[urr]_{j_0} \ar[ull]^{\pi^E} &&}
\end{equation}

We define the following kernels
\[
\widetilde{S}_1 = \widetilde{i}_+ \mco_{\widetilde{V}}, \quad \widetilde{S}_Y = \widetilde{i}_+j_{\dot{V} \dag} \mco_{\dot{V}}, \quad \widetilde{S}_{1/t} = \widetilde{i}_+j_{\dot{V} +} \mco_{\dot{V}}, \quad \widetilde{S}_\delta = \widetilde{i}_+ i_{E +} \mco_E
\]
on $\mbp \times V$.

The following result, comparing the affine Radon transform and the Fourier-Laplace transform, is proven in \cite{AE} for the case $k=1$.
\begin{prop}\cite[Proposition 1]{AE}\label{prop:isoAE} There is the following isomorphism of integral kernels in $D^b(\mcd_{\mbp \times \hat{V}})$.
\[
\widetilde{S}_{u} \diamond L \simeq  R'_{\hat{u}}\, ,
\]
where
\[
u = 1, \delta , Y , 1/t \quad \text{and} \quad \hat{u} = \delta, 1, 1/t, Y \, .
\]
\end{prop}

Define the following functors from $D^b_{qc}(\mcd_{\mbp})$ to $D^b_{qc}(\mcd_V)$:
\begin{align}
ext_\delta(M) &:= j_{0 +} (\pi^E)^+ M \simeq \widetilde{j}_+ i_{E +} i_E^+ \widetilde{\pi}^+ M \notag \\
ext_1(M) &:= \widetilde{j}_+ \widetilde{\pi}^+ M \notag \\
ext_{1/t}(M) &:= j_+ \pi^+ M \simeq \widetilde{j}_+ j_{\dot{V} +} j_{\dot{V}}^+ \widetilde{\pi}^+(M) \notag
\end{align}
resp. from $D^b_{coh}(\mcd_\mbp)$ to $D^b_{coh}(\mcd_V)$
\[
ext_{Y}(M) := j_\dag \pi^+ M \simeq \widetilde{j}_+ j_{\dot{V} \dag} j_{\dot{V}}^+ \widetilde{\pi}^+(M)
\]

\begin{lem}
Let $M \in  D^b_{coh}(\mcd_{\mbp})$, we have the following triangles in $D^b_{coh}(\mcd_V)$:
\begin{enumerate}
\item $ext_1(M) \lra ext_{1/t}(M) \lra ext_{\delta}(M) \overset{+1}{\lra}$
\item $ext_{\delta}(M) \lra ext_Y(M) \lra ext_{1}(M) \overset{+1}{\lra}$
\end{enumerate}
\end{lem}
\begin{proof}
The triangles can be deduced from
\[
i_{E +} i_{E}^+[-1] (M) \lra  M \lra j_{\dot{V} + } j_{\dot{V}}^+(M) \overset{+1}{\lra}
\]
and the corresponding dual triangle.
\end{proof}
The proposition above gives the following comparison between the affine Radon transform and the Fourier-Laplace transform.

\begin{cor}\label{cor:affRad}
For $M \in D^b_{qc}(\mcd_\mbp)$ there are natural isomorphisms in $D^b_{qc}(\mcd_{\hat{V}}) $
\[
\FL \circ ext_u (M) \simeq \mcr'_{\hat{u}}(M)\, \quad u = 1, \delta
\]
and for $M \in D^b_{coh}(\mcd_\mbp)$ there are natural isomorphisms in $D^b_{coh}(\mcd_{\hat{V}})$
\[
\FL \circ ext_u (M) \simeq \mcr'_{\hat{u}}(M)\, \quad u = 1/t, Y\, .
\]

Furthermore, there are the following isomorphisms of triangles
\[
\xymatrix{\mcr'_1(M) \ar[r]  \ar[d]^\simeq & \mcr'_{1/t}(M) \ar[r] \ar[d]^\simeq  & \mcr'_{\delta}(M) \ar[r]^{+1} \ar[d]^\simeq  & \\
FL \circ ext_\delta(M) \ar[r] & FL \circ ext_{Y}(M) \ar[r] & FL \circ ext_1(M) \ar[r]^<<<<{+1} & }
\]
and
\[
\xymatrix{\mcr'_\delta(M) \ar[r]  \ar[d]^\simeq & \mcr'_{Y}(M) \ar[r] \ar[d]^\simeq  & \mcr'_{1}(M) \ar[r]^{+1} \ar[d]^\simeq  & \\
FL \circ ext_1(M) \ar[r] & FL \circ ext_{1/t}(M) \ar[r] & FL \circ ext_\delta(M) \ar[r]^<<<<{+1} & }
\]
\end{cor} 
\begin{proof}
We have
\[
\mcr'_{\hat{u}}(M) \simeq M \diamond R'_{\hat{u}} \simeq M \diamond \widetilde{S}_u \diamond L \simeq FL(M \diamond \widetilde{S}_u)\, .
\]
So we need to shows that $ext_u(M) \simeq M \diamond \widetilde{S}_u$. But this follows from \cite[Lemma 1]{AE}. The isomorphism of triangles follows from the fact that the proof of \cite[Proposition 1]{AE} is completely functorial.
\end{proof}

The two results above which deal with the affine Radon transform were used by d'Agnolo and Eastwood to prove the following theorem, which generalizes a result obtained by Brylinski in \cite[Th\'{e}or\`{e}me 7.27]{Brylinski}

\begin{thm}\cite[Theorem 2]{AE}\label{thm:AEthm}
For $M \in D^b_{qc}(\mcd_\mbp)$ there are natural isomorphisms in $D^b_{qc}(\mcd_{\hat{V} \setminus \{ 0 \}}) $
\[
r^+ \FL \circ ext_u (M) \simeq \hat{\pi}^+ \mcr_{\hat{u}}(M) \quad u= 1,\delta
\]
and for $M \in D^b_{coh}(\mcd_\mbp)$ there are natural isomorphisms in $D^b_{coh}(\mcd_{\hat{V} \setminus \{ 0 \}}) $
\[
r^+ \FL \circ ext_u (M) \simeq \hat{\pi}^+ \mcr_{\hat{u}}(M) \quad u = 1/t, Y\, .
\]
where $r: \hat{V} \setminus \{0\} \ra \hat{V}$ is the natural inclusion and $\hat{\pi}: \hat{V} \setminus \{0\} \ra  \mbp$ the canonical projection.
\end{thm}

\subsection{Statement of results}\label{sec:statresult}
In order to state the comparison theorem in the case $k >1$ we have to introduce the following maps. Let
\[
\Delta : \mbp \lra \mbp^{\times k}
\]
be the diagonal embedding. Consider the following cartesian diagram
\[
\xymatrix{\dot{\mbv} \ar[r] \ar[d]_{\pi_k} & \dot{V}^{\times k} \ar[d]^{\prod \pi}\\ \mbp \ar[r]^\Delta & \mbp^{\times k}}
\]
where the vertical map $\prod \pi$ on the right is just the $k$-th product map of the canonical projection $\pi : \dot{V} \ra \mbp$ and $\dot{\mbv}$ is defined as the $k$-times product:
\begin{equation}\label{eq:foldedvdot}
\dot{\mbv} := \dot{V} \times_\mbp \ldots \times_\mbp \dot{V} \, .
\end{equation}
Recall that the closure of the map $(\pi,j): \dot{V} \ra \mbp \times V$ is the space $\widetilde{V}$ defined above which is the total space $\mbv(\mco_\mbp(-1))$ of the tautological line bundle $\mco_{\mbp}(-1)$. We get the following commutative diagram
\[
\xymatrix{\mbp \times V^{\times k} \ar[r] & (\mbp \times V)^{\times k} \\ \dot{\mbv} \ar[u] \ar[r] & \dot{V}^{\times k } \ar[u]_{\prod(\pi,j)} }
\]
where the closure of the image of $\prod(\pi,j)$ is equal to $\widetilde{V}^{\times k}$. Hence the closure $\mbv$ of the image of $\dot{\mbv} \ra \mbp \times V^{\times k}$  is equal to 
\[
\mbv = \mbv\left(\bigoplus_{i=1}^k \mco_\mbp(-1)\right)
\]
which is the total space of the $k$-times direct sum of the tautological bundle of $\mbp$. Denote by $D$ the complement of $\dot{\mbv}$ in $\mbv $, by $E$ the image of the zero section in $\mbv$ and by $\mbv^\circ$ the complement of $E$ in $\mbv$. We get the following diagrams
\[
\xymatrix{ && \mbv^\circ \ar[d]^{j_{\mbv^\circ}} \ar[drr]^{j^\circ} \ar[dll]_{\pi^\circ} && \\
\mbp && \mbv \ar[rr]^{j_k} \ar[ll]_{\pi_k} && V^{\times k} \\
&&  E \ar[u]_{i_{E}} \ar[urr]_{j_E} \ar[ull]^{\pi^E} &&} \qquad \xymatrix{ && \dot{\mbv} \ar[d]^{j_{\dot{\mbv}}} \ar[drr]^{j^\ast} \ar[dll]_{\pi^\ast} && \\
\mbp && \mbv \ar[rr]^{j_k} \ar[ll]_{\pi_k} && V^{\times k}\, , \\
&&  D \ar[u]_{i_{D}} \ar[urr]_{j_{D}} \ar[ull]^{\pi^D} &&}
\]
where $\pi^\circ, \pi^E, \pi^\ast$ and  $\pi^D$ are  restrictions of the projection $\mbp \times V^{\times k} \ra \mbp$ to the first factor and $j^\circ, j^\ast, j_E$ and $j_D$ are restrictions of the projection to the second factor.

Define the following functors from $D^b_{qc}(\mcd_{\mbp})$ to $D^b_{qc}(\mcd_{V^{\times k}})$:
\begin{align}
ext^k_\delta(M) &:= (j_{E})_+ (\pi^E)^+ M \simeq (j_{k})_+ (i_{E})_+ (i_{E})^+ (\pi_k)^+ M \notag \\
ext_1^k(M) &:= (j_{k})_+ (\pi_k)^+ M \notag \\
ext^k_{1/t}(M) &:= (j^\circ)_+ (\pi^\circ)^+ M \simeq (j_{k})_+ (j_{\mbv^\circ})_+ (j_{\mbv^\circ})^+ (\pi_k)^+ M \notag \\
%ext^k_{C^+}(M) &:= (j^*_{k})_+ (\pi^*_k)^+ M \simeq (j_{k})_+ (j_{\dot{V},k})_+ (j_{\dot{V},k})^+ (\widetilde{\pi}_k)^+ M \notag 
\end{align}
resp. from $D^b_{coh}(\mcd_{\mbp})$ to $D^b_{coh}(\mcd_{V^{\times k}})$

\begin{align}
ext^k_{Y}(M) &:= (j^\circ)_\dag (\pi^\circ)^+ M \simeq (j_k)_+ (j_{\mbv^\circ})_\dag (j_{\mbv^\circ})^+ (\pi_k)^+ M \notag \\
%ext^k_{A^\dag}(M) &:=(j_{k})_\dag\, (\mbd \circ  R\Gamma_{E^k} \circ \mbd) \, (\pi_k)^+ M \notag \\
ext^k_{U^\dag}(M) &:= (j^*)_\dag (\pi^*)^+ M \simeq (j_{k})_+ (j_{\dot{\mbv}})_\dag (j_{\dot{\mbv}})^+ (\widetilde{\pi}_k)^+ M \notag 
\end{align}

We still need one more extension functor. Let $\mbv' := \mbv\left( \bigoplus_{i=1}^k \mco_\mbp(1)\right)$ be the dual vector bundle of $\mbv$.
Denote by $D'_l$ for $l=1, \ldots ,k$ the image of the zero section of the $l$-the summand $\mco_{\mbp}(1)$ and by $D' = \bigcup_{i=1}^l D'_l$ their union in $\mbv'$. Let $\dot{\mbv}' \overset{j_{\dot{\mbv}'}}{\lra} \mbv'$ the complement of $D'$. We get the following triangle in $D^b_{coh}(\mcd_{\mbv'})$
\[
R\Gamma_{D'}\mco_{\mbv'} \lra \mco_{\mbv'} \lra (j_{\dot{\mbv}'})_+ (j_{\dot{\mbv}'})^+ \mco_{\mbv'} \overset{+1}{\lra}
\]
We now apply the Fourier-Laplace transformation $FL_\mbp$ to the triangle above:
\begin{equation}\label{eq:triFL}
FL_\mbp(R\Gamma_{D'}\mco_{\mbv'}) \lra i_{E +} \mco_{E} \lra (j_{\dot{\mbv}})_\dag (j_{\dot{\mbv}'})^+ \mco_{\mbv'} \overset{+1}{\lra}
\end{equation}

In order to see this we have to compute the second and third term. We will do this locally with respect to $\mbp$. Let $(v_0: \ldots : v_n)$ be homgeneous coordinates on $\mbp$ and let $V_i \subset \mbp$ the affine chart with coordinates $(v_{i1}, \ldots , v_{in})$ given by $v_i \neq 0$. Then $\mbv_i := \mbv_{\mid V_i} \simeq V_i \times \mbc^k$ and $ \mbv'_i := \mbv'_{\mid V_i} \simeq V_i \times (\mbc^k)'$ are trivial vector bundles with fiber coordinates $\xi_1, \ldots , \xi_k$  and $x_1, \ldots , x_k$, respectively.  Hence by Proposition \ref{prop:FLconc} we have the following isomorphisms in $D^b_{qc}(\mcd_{\mbv_i})$:
\begin{align}
&FL_{V_i}(\mco_{\mbv_i'} ) \simeq FL_{V_i}\left(\mcd_{\mbv_i'}/(\p_{v_{i1}}, \ldots \p_{v_{i_n}},\p_{x_1}, \ldots ,\p_{x_{k}})\right)\notag \\ \simeq  &\; \mcd_{\mbv_i} / \left(\p_{v_{i1}}, \ldots , \p_{v_{in}}, \xi_1, \ldots \xi_k) \right) \simeq i_{E_i + } \mco_{E_i}\, , \notag
\end{align}
where $i_{E_i}: E_i := E \cap \mbv_i \ra \mbv_i $ is the canonical inclusion of the zero section.
%The second term follows from
%\[
%FL(\mco_{\mbc^n}) = i_+ \mco_{\{0\}}
%\]
%where $i: \{0\} \ra \mbc^n$ is the canonical inclusion and 

Denote by $j_{\dot{\mbv}'_i}: \dot{\mbv}'_i := \dot{\mbv}' \cap \mbv'_i \ra \mbv'_i$ and $j_{\dot{\mbv}_i}: \dot{\mbv}_i := \dot{\mbv} \cap \mbv_i \ra \mbv_i$ the canonical inclusions. The third term follows from the following isomorphisms in $D^b_{qc}(\mcd_{\mbv_i})$:
\begin{align}
&FL_{V_i}\left((j_{\dot{\mbv}'_i})_+ (j_{\dot{\mbv}'_i})^+ \mco_{\mbv'_i}\right) \simeq FL_{V_i}\left(\mco_{\mbv'_i}(\ast (D' \cap \mbv'_i))\right) \simeq FL_{V_i}(\mcd_{\mbv'_i}/ (\p_{v_{i1}}, \ldots , \p_{v_{in}},\p_{x_1}x_1, \ldots , \p_{x_n} x_n)) \notag \\
\simeq\; &\mcd_{\mbv_i}/ (\p_{v_{i1}}, \ldots , \p_{v_{in}},\xi_1\p_{\xi_1}, \ldots , \xi_n\p_{\xi_n} ) \simeq (j_{\dot{\mbv}_i})_\dag (j_{\dot{\mbv}_i})^+ \mco_{\mbv_i}\, , \notag
\end{align}
where the first isomorphism follows from $(j_{\dot{\mbv}'_i})_+ (j_{\dot{\mbv}'_i})^+ \mco_{\mbv'_i} \simeq (j_{\dot{\mbv}'_i})_* (j_{\dot{\mbv}'_i})^{-1}\mco_{\mbv'_i}$ (since $j_{\dot{\mbv}'_i}$ is an open embedding) and the last isomorphism follows from
\begin{align}
&(j_{\dot{\mbv}_i})_\dag (j_{\dot{\mbv}_i})^+ \mco_{\mbv_i} \simeq \mbd (j_{\dot{\mbv}_i})_+ \mbd (j_{\dot{\mbv}_i})^{-1} \mco_{\mbv_i} \simeq \mbd\, (j_{\dot{\mbv}_i})_*\, \mbd\, \mco_{\dot{\mbv}_i} \simeq  \mbd\, (j_{\dot{\mbv}_i})_*\, \mco_{\dot{\mbv}_i} \simeq \mbd\, \mco_{\mbv_i}(\ast D \cap \mbv_i) \notag \\
\simeq\; &\mbd\,\left( \mcd_{\mbv_i}/ (\p_{v_{i1}}, \ldots , \p_{v_{in}},\p_{\xi_1}\xi_1, \ldots , \p_{\xi_n}\xi_n ) \right) \simeq \mcd_{\mbv_i}/ (\p_{v_{i1}}, \ldots , \p_{v_{in}},\xi_1\p_{\xi_1}, \ldots , \xi_n\p_{\xi_n} )\, . \notag
\end{align}
%\[
%FL(\mch^0(j_+ \mco_{(\mbc^*)^n})) \simeq FL(\mcd_{\mbc^n} / (\p_i \lambda_i)_{i=1,\ldots ,n}) = FL(\mcd_{\mbc^n} / (\lambda_i\p_i)_{i=1,\ldots ,n}) \simeq FL(\mch^0(j_\dag \mco_{(\mbc^*)^n}))
%\]
We are now able to define the last extension functor:
\[
ext^k_{A^\dag}(M) :=(j_{k})_+\, ( (\pi_k)^+ M \overset{L}{\otimes}FL_\mbp(R\Gamma_{D'}\mco_{\mbv'}))
\]
from $D^b_{coh}(\mcd_\mbp)$ to $D^b_{coh}(\mcd_{V^{\times k}})$.
\begin{lem}
Let $M \in D^b_{coh}(\mcd_\mbp)$, whe have the following triangles in $D^b_{coh}(\mcd_{V^{\times k}})$:
\begin{enumerate}
\item $ext^k_1(M) \lra ext^k_{1/t}(M) \lra ext^k_{\delta}(M) \overset{+1}{\lra}$
\item $ext^k_{\delta}(M) \lra ext^k_Y(M) \lra ext^k_{1}(M) \overset{+1}{\lra}$
\item $ext^k_{A^\dag}(M) \lra ext^k_\delta(M) \lra ext^k_{U^\dag}(M) \overset{+1}{\lra}$
\end{enumerate}
\end{lem}
\begin{proof}
The first triangle can be deduced from
\[
i_{E +} i_{E}^+[-1] (M) \lra  M \lra j_{\dot{\mbv} + } j_{\dot{\mbv}}^+(M) \overset{+1}{\lra}
\]
and the second from corresponding dual triangle. The third triangle follows from \eqref{eq:triFL}.
\end{proof}
% \notag \\
%
%where $N_A := FL_\mbp
%\renewcommand{\arraystretch}{1.4} 
%\begin{tabular}{|c|c|c|c|c|c|c|}
%\hline
%$u$ & $1$ & $\delta$ & $1/t$ & $Y$ & $U^\dag$ & $A^\dag$ \\
%\hline
%$\hat{u}$ &$\delta$ & $1$ & $Y$ & $1/t$ & $U^+$ & $A^+$\\
%\hline
%\end{tabular}

We have the following generalization of Corollary \ref{cor:affRad}
\begin{prop}\label{prop:higheraffRad}
For $M \in D^b_{qc}(\mcd_{\mbp})$ there are natural isomorphisms in $D^b_{qc}(\mcd_{\hat{W}}) $
\[
\FL \circ ext_u^k (M) \simeq \mcr'_{\hat{u}}(M)\, \quad u = 1, \delta \quad \hat{u} =\delta ,1
\]
and for $M \in D^b_{coh}(\mcd_\mbp)$ there are natural isomorphisms in $D^b_{coh}(\mcd_{\hat{W}})$
\[
\FL \circ ext^k_u (M) \simeq \mcr'_{\hat{u}}(M)\, \quad u = 1/t, Y, U^\dag, A^\dag \quad \hat{u} = Y, 1/t, U^+, A^+ \, .
\]
Furthermore, there are the following isomorphisms of triangles
\[
\xymatrix{\mcr'_1(M) \ar[r]  \ar[d]^\simeq & \mcr'_{1/t}(M) \ar[r] \ar[d]^\simeq  & \mcr'_{\delta}(M) \ar[r]^{+1} \ar[d]^\simeq  & \\
FL \circ ext^k_\delta(M) \ar[r] & FL \circ ext^k_{Y}(M) \ar[r] & FL \circ ext^k_1(M) \ar[r]^<<<<{+1} & }
\]
\[
\xymatrix{\mcr'_\delta(M) \ar[r]  \ar[d]^\simeq & \mcr'_{Y}(M) \ar[r] \ar[d]^\simeq  & \mcr'_{1}(M) \ar[r]^{+1} \ar[d]^\simeq  & \\
FL \circ ext_1(M) \ar[r] & FL \circ ext_{1/t}(M) \ar[r] & FL \circ ext_\delta(M) \ar[r]^<<<<{+1} & }
\]
\[
\xymatrix{\mcr'_{A^+}(M) \ar[r]  \ar[d]^\simeq & \mcr'_{1}(M) \ar[r] \ar[d]^\simeq  & \mcr'_{U^+}(M) \ar[r]^{+1} \ar[d]^\simeq  & \\
FL \circ ext_{A^\dag}(M) \ar[r] & FL \circ ext_{\delta}(M) \ar[r] & FL \circ ext_{U^\dag}(M) \ar[r]^<<<<{+1} & }
\]

\end{prop}

The corresponding generalization  of Theorem \ref{thm:AEthm} is the following theorem.
\begin{thm}\label{thm:higherRad}
For $M \in D^b_{qc}(\mcd_{\mbp})$ there are natural isomorphisms in $D^b_{qc}(\mcd_{S(k,n)})$
\[
r^+ (\FL \circ ext^k_u (M)) \simeq \hat{\pi}^+ \mcr_{\hat{u}}(M) \quad u= 1,\delta\quad \hat{u} = \delta, 1
\]
and for $M \in D^b_{coh}(\mcd_{\mbp})$ there are natural isomorphisms in $D^b_{coh}(\mcd_{S(k,n)}) $
\[
r^+ (\FL \circ ext^k_u (M)) \simeq \hat{\pi}^+ \mcr_{\hat{u}}(M) \quad u = 1/t, Y \quad \hat{u} = Y, 1/t \, .
\]
where $r: S(k,n) \ra \hat{W} = \hat{V}^{\times k}$ is the natural inclusion and $\hat{\pi}: S(k,n) \ra  \mbg$ the canonical projection.
\end{thm}

\subsection{Proof of the higher codimension case}

In this section we give the proof of the comparison theorem in the case of $k > 1$, using the case $k=1$ already proven by d'Agnolo and Eastwood.

We first introduce various Radon like transformations on the product space $\mbp^{\times k}$ and prove in Lemma \ref{lem:RadontoRadonprod} that these transformations applied to some $\mcd$-module $M$ which is supported on the diagonal of $\mbp^{\times k}$ are isomorphic to the $d$-plane Radon transforms of $M$. Set
\[
\widetilde{Z}:= \{ ([v^1], \ldots ,[v^k] , \lambda^1, \ldots , \lambda^k) \in \mbp^{\times k} \times \hat{W} \mid \lambda^i(v^i) = 0\ \forall i \in \{1,\ldots ,k\}\}
\]
and denote by $i_{\widetilde{Z}}: \widetilde{Z} \lra \mbp^{\times k} \times \hat{W}$ its inclusion. 
Define the open subvariety 
\[
\widetilde{U} := \{ ([v^1], \ldots , [v^k], \lambda^1, \ldots , \lambda^k \in \mbp^{\times k} \times \hat{W} \mid \lambda^1(v^1) \neq 0, \ldots , \lambda^k(v^k) \neq 0  \}
\]
and let $j_{\widetilde{U}}: \widetilde{U} \lra \mbp^{\times k} \times \hat{W}$ be the canonical inclusion.
Denote by $\pi_1^\Delta, \pi_2^\Delta$ the projections from $\mbp^{\times k} \times \hat{W}$ to $\mbp^{\times k}$ resp. $\hat{W}$ 
and denote their restrictions to $\widetilde{Z}$ resp. $\widetilde{U}$ by
$\pi_1^{\widetilde{Z}}, \pi_2^{\widetilde{Z}}$ resp. $\pi_1^{\widetilde{U}}, \pi_2^{\widetilde{U}}$.

We define the following transformations for $M \in D^b_{qc}(\mcd_{\mbp^{\times k}})$
\begin{align}
\widetilde{\mcr}_{\delta}(M) &:= \pi_{2+}^{\widetilde{Z}} (\pi_1^{\widetilde{Z}})^+ M \simeq \pi_{2+}^{\Delta} i_{\widetilde{Z}+} i_{\widetilde{Z}}^+ (\pi_1^{\Delta})^+ M\, , \notag \\
\widetilde{\mcr}_{U^+}(M)&:= \pi_{2+}^{\widetilde{U}} (\pi_1^{\widetilde{U}})^+ M \simeq \pi_{2+}^{\Delta} j_{\widetilde{U}+} j_{\widetilde{U}}^+ (\pi_1^{\Delta})^+ M\, , \notag \\
\widetilde{\mcr}_1(M) &:= \pi_{2+}^{\Delta} (\pi_1^\Delta)^+ M\, . \notag
\end{align}
%and for $M \in D^b_{coh}(\mcd_{\mbp^{\times k}})$
%\[
%\widetilde{\mcr}_{U^\dag}(M) := \pi_{2\dag}^{\widetilde{U}} (\pi_1^{\widetilde{U}})^+ M \simeq \pi_{2+}^{\Delta} j_{\widetilde{U}\dag} j_{\widetilde{U}}^+ (\pi_1^{\Delta})^+ M\, .
%\]
Consider the following diagram
$$
\xymatrix{&&& \hat{W} &\\ U'  \ar[rr]_{j_{U'}} \ar[dd]_{\Delta_{U'}}  && \mbp \times \hat{W} \ar[ur]^{\pi_2} \ar[dl]_{\pi_1} \ar[dd]_{\Delta \times id_{\hat{W}}} && Z'  \ar[dd]_{\Delta_{Z'}} \ar[ll]^/-1em/{i_{Z'}}  \\ &\mbp \ar[dd]_/-1em/{\Delta} & && \\ \widetilde{U} \ar[rr]_/.5em/{j_{\widetilde{U}}}|<<<<<<<<<<<<<\hole\hole && \mbp^{\times k} \times \hat{W} \ar[uuur]_{\pi_2^\Delta}|<<<<<<<<<<<<<<<<<<<<<<<<<<<\hole \ar[dl]^{\pi_1^\Delta}&& \widetilde{Z} \ar[ll]_{i_{\widetilde{Z}}} \\ &  \mbp^{\times k} &&& }
$$
where all squares are Cartesian.
\begin{lem}\label{lem:RadontoRadonprod}
We have the following isomorphism for $M \in D^b_{qc}(\mcd_{\mbp})$
\[
\mcr'_{\hat{u}}(M) \simeq \widetilde{\mcr}_{\hat{u}}(\Delta_+(M)) \quad \hat{u} = \delta,1,U^+\, .
\]
%and for $M \in D^b_{coh}(\mcd_{\mbp})$
%\[
%\mcr'_{U^\dag}(M) \simeq \widetilde{\mcr}_{U^\dag}(\Delta_+(M))\, .
%\]
\end{lem}
\begin{proof}
We are going to show $\mcr'_\delta(M) \simeq \widetilde{\mcr}_\delta(\Delta_+(M))$ the other cases are similar or even simpler.
\begin{align}
\mcr'_{\delta}(M) &\simeq \pi_{2+} i_{Z' +} i_{Z'}^+ \pi_1^+ M \notag  \\
&\simeq \pi_{2 +}^{\Delta} (\Delta \times id_w)_+ i_{Z' +} i_{Z'}^+ \pi_1^+ M \notag  \\
&\simeq \pi_{2 +}^{\Delta}i_{\widetilde{Z} +} \Delta_{Z' +} i_{Z'}^+ \pi_1^+ M \notag  \\
&\simeq \pi_{2 +}^{\Delta} i_{\widetilde{Z} +} i_{\widetilde{Z}}^+ (\Delta \times id_W)_+ \pi_1^+ M \notag  \\
&\simeq \pi_{2 +}^{\Delta} i_{\widetilde{Z} +} i_{\widetilde{Z}}^+ \pi_{1+}^\Delta \Delta_+ M \notag  \\
&= \widetilde{\mcr}_\delta(\Delta_+(M))\, . \notag
\end{align}
\end{proof}
Now we want to define transformations $\widetilde{\mcr}^i_{\hat{u}}$ which can be considered as partial Radon transforms with respect to the $i$-th factor of $\mbp^{\times k}$. We will prove in Lemma \ref{lem:RadonprodtoRadonpart} that the Radon like transformations $\widetilde{R}_{\hat{u}}$ introduced above are actually equivalent to the consecutive application of the corresponding hyperplane Radon transforms $\widetilde{\mcr}^i_{\hat{u}}$ on each factor. 

Define $T_{i,j}:= \mbp^{\times i} \times \hat{V}^{\times j}$ for $i,j \in \{0, \ldots , k\}$ with $k \leq i+j \leq 2 \cdot k$ and consider the following diagram of spaces
$$
\xymatrix@R=2.5em{T_{k,k} \ar@{.}[rrrr] & & & & T_{1,k} \ar[d]_{\pi_1^{T_{1,k}}} \ar[r]_{\pi_2^{T_{1,k}}} & T_{0,k} \\ & & & &  T_{1,k-1} &\\ & & & & & \\ T_{k,2} \ar[r]^{\pi_2^{T_{k,2}}} \ar[d]_{\pi_1^{T_{k,2}}} \ar@{.}[uuu]& T_{k-1,2} \ar[r]^{\pi_2^{T_{k-1,2}}} \ar[d]_{\pi_1^{T_{k-1,2}}}& T_{k-2,2} \ar@{.}[uurr]\\ T_{k,1} \ar[r]_{\pi_2^{T_{k,1}}} \ar[d]_{\pi_1^{T_{k,1}}} & T_{k-1,1} \\ T_{k,0} }
$$

where the maps are defined as follows:
\begin{align}
\pi_1^{T_{i,j}}: T_{i,j} &\lra T_{i,j-1}\, , \notag \\
([v^1], \ldots , [v^i], \lambda^1, \ldots , \lambda^j) &\mapsto ([v^1], \ldots , [v^i], \lambda^2, \ldots , \lambda^j)\, ,\notag \\
\pi_2^{T_{i,j}}: T_{i,j} &\lra T_{i-1,j}\, , \notag \\
([v^1], \ldots , [v^i], \lambda^1, \ldots , \lambda^j) &\mapsto ([v^1], \ldots , [v^{i-1}], \lambda^1, \ldots , \lambda^j)\notag
\end{align}
for $i,j \in \{ 1, \ldots ,k\}$ with $k+1 \leq i+j \leq 2 \cdot k$.
Now it is easy to see that the squares
$$
\xymatrix{ T_{i,j} \ar[r]^{\pi_2^{T_{i,j}}} \ar[d]_{\pi_1^{T_{i,j}}} & T_{i-1,j} \ar[d]^{\pi_1^{T_{i-1,j}}}\\ T_{i,j-1} \ar[r]^{\pi_2^{T_{i,j-1}}} & T_{i-1,j-1}}
$$
are cartesian for $i,j \in\{2,\ldots,k\}$ with $k+2 \leq i+j \leq 2 \cdot k$.

\noindent We now define closed subvarieties $Z_{i,j} \subset T_{i,j}= \mbp^{\times i} \times \hat{V}^{\times j}$ for $i,j \in \{0, \ldots ,k\}$ with $k+1 \leq i+j \leq 2\cdot k$:
\[
Z_{i,j} := \{([v^1],\ldots , [v^i], \lambda^1, \ldots , \lambda^j) \in T_{i,j} \mid \lambda^1(v^{k-j+1}) = \ldots = \lambda^{i+j-k}(v^{i}) =0\}\,.
\]
Notice that 
\[
Z_{k,k} = \{ [v^1], \ldots , [v^k], \lambda^1, \ldots , \lambda^k \in T_{k,k} \mid \lambda^1(v^1) = \ldots \lambda^k(v^k) = 0\} = \widetilde{Z}
\]
and
\[
Z_{i,k-i+1} = \{[v^1], \ldots , [v^i], \lambda^1, \ldots , \lambda^{k-i+1} \in T_{i, k-i+1} \mid \lambda^1(v^i)=0 \}
\]
for $i \in \{1, \ldots , k\}$.

\begin{lem}
Let $\pi_1^{Z_{i,j}}$ resp. $\pi_2^{Z_{i,j}}$ be the restrictions of $\pi_1^{T_{i,j}}$ resp. $\pi_2^{T_{i,j}}$ to $Z_{i,j}$. The squares
$$
\xymatrix{ Z_{i,j} \ar[dd]_{\pi_1^{Z_{i,j}}} \ar[r]^{\pi_2^{Z_{i,j}}} & Z_{i-1,j} \ar[dd]^{\pi_1^{Z_{i-1,j}}} & & & Z_{r+1, s+1}\ar[r]^{\pi_2^{Z_{r+1, s+1}}} \ar[dd]_{\pi_1^{Z_{r+1, s+1}}}&Z_{r,s+1}\ar[dd]^{\pi_1^{Z_{r,s+1}}}\\
& &  \text{and} & & & \\
Z_{i, j-1} \ar[r]^{\pi_2^{Z_{i,j-1}}} & Z_{i-1,j-1} &  & & Z_{r+1,s}\ar[r]^{\pi_2^{Z_{r+1,s}}}&T_{r,s}}
$$
are cartesian for $i,j \in \{2, \ldots ,k\}$ with $k+3 \leq i+j \leq 2 \cdot k$ and for $r,s \in \{1, \ldots ,k-1\}$ with $r+s = k$. 
\end{lem}
\begin{proof}
First, we write down the definition of the spaces involved
\begin{align}
Z_{i,j} &= \{([v^1],\ldots , [v^i], \lambda^1, \ldots , \lambda^j) \in T_{i,j} \mid \lambda^1(v^{k-j+1}) = \ldots = \lambda^{i+j-k}(v^{i}) =0\}\, , \notag \\
Z_{i-1,j} &= \{([v^1],\ldots , [v^{i-1}], \lambda^1, \ldots , \lambda^j) \in T_{i-1,j} \mid \lambda^1(v^{k-j+1}) = \ldots = \lambda^{i-1+j-k}(v^{i-1}) =0\}\, , \notag \\
Z_{i,j-1} &= \{([v^1],\ldots , [v^i], \lambda^2, \ldots , \lambda^j) \in T_{i,j-1} \mid \lambda^2(v^{k-j+2}) = \ldots = \lambda^{i+j-k}(v^{i}) =0\}\, , \notag \\
Z_{i-1,j-1} &= \{([v^1],\ldots , [v^{i-1}], \lambda^2, \ldots , \lambda^{j-1}) \in T_{i-1,j-1} \mid \lambda^2(v^{k-j+2}) = \ldots = \lambda^{i-1+j-k}(v^{i-1}) =0\}\, , \notag
\end{align}
where we changed the indices of the elements in the definition of the spaces $Z_{i,j-1}$ and $Z_{i-1,j-1}$.
Now one can see rather easily that the maps
$$
\xymatrix{Z_{i,j} \ar[d]_{\pi_1^{Z_{i,j}}} \ar[r]^{\pi_2^{Z_{i,j}}} & Z_{i-1,j} \\ Z_{i, j-1} &  }
$$
are well-defined and that the left square is cartesian. The proof of the fact that the right square is cartesian is completely parallel if one sets $Z_{i,j} := T_{i,j}$ for $i+j =k$.
\end{proof}

The lemma above gives rise to the following diagram of cartesian squares
\begin{equation}\label{Zdiag}
\xymatrix{Z_{k,k} \ar@{.}[rrrr] & & & & Z_{1,k} \ar[d]_{\pi_1^{Z_{1,k}}} \ar[r]_{\pi_2^{Z_{1,k}}} & T_{0,k} \\ & & & &  T_{1,k-1} &\\ & & & & & \\ Z_{k,2} \ar[r]^{\pi_2^{Z_{k,2}}} \ar[d]_{\pi_1^{Z_{k,2}}} \ar@{.}[uuu]& Z_{k-1,2} \ar[r]^{\pi_2^{Z_{k-1,2}}} \ar[d]_{\pi_1^{Z_{k-1,2}}}& T_{k-2,2} \ar@{.}[uurr]\\ Z_{k,1} \ar[r]_{\pi_2^{Z_{k,1}}} \ar[d]_{\pi_1^{Z_{k,1}}} & T_{k-1,1} \\ T_{k,0} }
\end{equation}

\noindent We now define open subvarieties $U_{i,j} \subset T_{i,j}= \mbp^{\times i} \times \hat{V}^{\times j}$ for $i,j \in \{0, \ldots ,k\}$ with $k+1 \leq i+j \leq 2\cdot k$:
\[
U_{i,j} := \{([v^1],\ldots , [v^i], \lambda^1, \ldots , \lambda^j) \in T_{i,j} \mid \lambda^1(v^{k-j+1})\neq 0,\ldots, \lambda^{i+j-k}(v^{i}) \neq 0\}\,.
\]
Notice that 
\[
U_{k,k} = \{ [v^1], \ldots , [v^k], \lambda^1, \ldots , \lambda^k \in T_{k,k} \mid \lambda^1(v^1)\neq 0,\ldots, \lambda^k(v^k) \neq 0\} =\widetilde{U}
\]
and
\[
U_{i,k-i+1} = \{[v^1], \ldots , [v^i], \lambda^1, \ldots , \lambda^{k-i+1} \in T_{i, k-i+1} \mid \lambda^1(v^i) \neq 0 \}
\]
for $i \in \{1, \ldots , k\}$. Arguing as in the case of the $Z_{i,j}$ we get a diagram of cartesian squares

$$
\xymatrix{U_{k,k} \ar@{.}[rrrr] & & & & U_{1,k} \ar[d]_{\pi_1^{U_{1,k}}} \ar[r]_{\pi_2^{U_{1,k}}} & T_{0,k} \\ & & & &  T_{1,k-1} &\\ & & & & & \\ U_{k,2} \ar[r]^{\pi_2^{U_{k,2}}} \ar[d]_{\pi_1^{U_{k,2}}} \ar@{.}[uuu]& U_{k-1,2} \ar[r]^{\pi_2^{U_{k-1,2}}} \ar[d]_{\pi_1^{U_{k-1,2}}}& T_{k-2,2} \ar@{.}[uurr]\\ U_{k,1} \ar[r]_{\pi_2^{U_{k,1}}} \ar[d]_{\pi_1^{U_{k,1}}} & T_{k-1,1} \\ T_{k,0} }
$$

We have the following diagram
\begin{equation}\label{eq:partRad}
\xymatrix{ && U_{i,j} \ar[drr]^{\pi_2^{U_{i,j}}} \ar[dll]_{\pi_1^{U_{i,j}}} \ar@{^(->}[d]^{j_{U_{i,j}}}&& \\
 T_{i,j-1} && T_{i,j} \ar[ll]_{\pi_1^{T_{i,j}}} \ar[rr]^{\pi_2^{T_{i,j}}} && T_{i-1,j} \\
  && Z_{i,j} \ar[ull]^{\pi_1^{Z_{i,j}}} \ar@{^(->}[u]_{i_{Z_{i,j}}} \ar[rru]_{\pi_2^{Z_{i,j}}} &&}
\end{equation}
for $i,j \in \{0, \ldots ,k\}$ with $i+j = k+1$. Now define the following partial Radon transforms for $M \in D^b_{qc}(\mcd_{T_{i,j-1}})$
\begin{align}
\widetilde{\mcr}^i_\delta(M)&:= \pi_{2 +}^{Z_{i,j}} (\pi_{1}^{Z_{i,j}})^+ M \simeq \pi_{2 +}^{T_{i,j}} i_{Z_{i,j} +} i_{Z_{i,j}}^+(\pi_{1}^{T_{i,j}})^+ M\, , \notag \\
\widetilde{\mcr}^i_{U^+}(M)&:= \pi_{2 +}^{U_{i,j}} (\pi_{1}^{U_{i,j}})^+ M \simeq \pi_{2 +}^{T_{i,j}} j_{U_{i,j} +} j_{U_{i,j}}^+(\pi_{1}^{T_{i,j}})^+ M\, , \notag \\
\widetilde{\mcr}^i_1(M)&:= \pi_{2 +}^{T_{i,j}} (\pi_{1}^{T_{i,j}})^+ M\, . \notag
\end{align}
%and for $M \in D^b_{coh}(\mcd_{T_{i,j-1}})$
%\[
%\widetilde{\mcr}^i_{Y}(M) := \pi_{2 \dag}^{U_{i,j}} (\pi_{1}^{U_{i,j}})^+ M \simeq \pi_{2 +}^{T_{i,j}} j_{U_{i,j} \dag} j_{U_{i,j}}^+(\pi_{1}^{T_{i,j}})^+ M\, .
%\]
\begin{lem}\label{lem:RadonprodtoRadonpart}
We have the following isomorphisms for $M \in D^b_{qc}(\mcd_{T_{k,0}})$
\[
\widetilde{\mcr}^1_{\hat{u}} \circ \ldots \circ \widetilde{\mcr}^k_{\hat{u}} (M) \simeq \widetilde{\mcr}_{\hat{u}}(M) \quad \hat{u} = 1, \delta, U^+ \, .
\]
%and for $M \in D^b_{coh}(\mcd_{T_{k,0}})$
%\[
%\widetilde{\mcr}^1_Y \circ \ldots \circ \widetilde{\mcr}^k_Y (M) \simeq \widetilde{\mcr}_Y(M)\, .
%\]
\end{lem}
\begin{proof}
We will only show $\widetilde{\mcr}^1_\delta \circ \ldots \circ \widetilde{\mcr}^k_\delta (M) \simeq \widetilde{\mcr}_\delta(M)$ the other cases are again similar or simpler. Recall the pyramid diagram \ref{Zdiag} from above.
The squares
$$
\xymatrix{ Z_{i,j} \ar[dd]_{\pi_1^{Z_{i,j}}} \ar[r]^{\pi_2^{Z_{i,j}}} & Z_{i-1,j} \ar[dd]^{\pi_1^{Z_{i-1,j}}} & & & Z_{r+1, s+1}\ar[r]^{\pi_2^{Z_{r+1, s+1}}} \ar[dd]_{\pi_1^{Z_{r+1, s+1}}}&Z_{r,s+1}\ar[dd]^{\pi_1^{Z_{r,s+1}}}\\
& &  \text{and} & & & \\
Z_{i, j-1} \ar[r]^{\pi_2^{Z_{i,j-1}}} & Z_{i-1,j-1} &  & & Z_{r+1,s}\ar[r]^{\pi_2^{Z_{r+1,s}}}&T_{r,s}}
$$
are cartesian for $i,j \in \{2, \ldots ,k\}$ with $k+3 \leq i+j \leq 2 \cdot k$ and for $r,s \in \{1, \ldots ,k-1\}$ with $r+s = k$. 

By base change we have the following isomorphisms of functors
\[
(\pi_{1}^{Z_{i-1,j}})^+ (\pi_{2}^{Z_{i,j-1}})_+ \simeq (\pi_{2}^{Z_{i,j}})_+ (\pi_{1}^{Z_{i,j}})^+ \quad \text{and} \quad (\pi_{1}^{Z_{r, s+1}})^+ (\pi_{2}^{Z_{r+1,s}})_+ \simeq (\pi_{2}^{Z_{r+1\, s+1}})_+ (\pi_{1}^{Z_{r+1 \, s+1}})^+\, .
\]
 Thus we have
\begin{align}
\widetilde{\mcr}^1_\delta \circ \ldots \circ \widetilde{\mcr}^k_\delta (M) &= (\pi_{2}^{Z_{1,k}})_+ (\pi_{1}^{Z_{1,k}})^+ \ldots (\pi_{2}^{Z_{k,1}})_+ (\pi_{1}^{Z_{k,1}})^+ M \notag \\
&\simeq (\pi_{2}^{Z_{1,k}})_+ \ldots (\pi_{2}^{Z_{k-1,k}})_+ (\pi_{2}^{Z_{k,k}})_+ (\pi_{1}^{Z_{k,k}})^+ \ldots (\pi_{1}^{Z_{k,2}})^+(\pi_{1}^{Z_{k,1}})^+ M\, . \notag
\end{align}
It now follows from the construction of the $Z_{i,j}$ that $Z_{k,k} \simeq \widetilde{Z}$ and, if we identify the latter spaces, that
\[
\pi_{2}^{Z_{1,k}}\ldots\pi_{2}^{Z_{k-1,k}}\pi_{2}^{Z_{k,k}} = \pi_2^{\widetilde{Z}} \quad \text{and} \quad \pi_{1}^{Z_{1,k}}\ldots\pi_{1}^{Z_{k, k-1}}\pi_{1}^{Z_{k,k}} = \pi_1^{\widetilde{Z}}\, .
\]
Thus
\begin{align}
\widetilde{\mcr}^1_\delta \circ \ldots \circ \widetilde{\mcr}^k_\delta (M) &\simeq \pi_{2+}^{\widetilde{Z}}(\pi_1^{\widetilde{Z}})^+ M \notag \\
&\simeq \widetilde{\mcr}_\delta (M)\, . \notag
\end{align}
\end{proof}
As a next step we want to compare the partial Radon transform $\widetilde{\mcr}^i_u$ with a partial Fourier-Laplace transform. For this, we need a description of the partial Radon transform as an integral transformation with kernel $\widetilde{R}^i_{\hat{u}}$. Recall diagram \eqref{eq:partRad} and set
\[
 \widetilde{R}^i_{U^+} := j_{U_{i,j} +} \mco_{U_{i,j}}, \quad \widetilde{R}^i_\delta := i_{Z_{i,j} +} \mco_{Z_{i,j}},\quad \widetilde{R}^i_1 := \mco_{T_{i,j}}\, . %\widetilde{R}^i_Y := j_{U_{i,j} \dag}\mco_{U_{i,j}}, \quad
\]
for $i \in \{1, \ldots , k\}$ and $j = k+1-i$. By arguing as in Proposition \ref{prop:kernelRadcompare}, we get
\[
M \diamond \widetilde{R}^i_{\hat{u}} := \pi_{2+}^{T_{i,j}}((\pi_1^{T_{i,j}})^{+}(M) \overset{L}{\otimes} \widetilde{R}^i_{\hat{u}}) \simeq \widetilde{\mcr}^i_u(M)
\]
for $M \in D^b_{qc}(\mcd_{T_{i,j-1}})$ and $\hat{u}= \delta,1,U^+$.\\% resp. $M \in D^b_{coh}(\mcd_{T_{i,j-1}})$ and $u=Y$.\\

Recall the definition of the space $T_{i,j}$ and define for $i+j = k+1$:
\[
T_i := T_{i,j} = \prod_{l=1}^{i-1} \mbp \times( \mbp \times \hat{V} )\times \prod_{l=i}^{k-1} \hat{V}
\] 
together with the following projection:
\begin{align}
\pi^{T_i}: T_i &\lra \mbp  \times \hat{V}\, ,  \notag \\
([v^1] ,\ldots [v^{i-1}], [v^i], \lambda^1, \lambda^{2}&,\ldots, \lambda^{k+1-i} ) \mapsto ([v^i], \lambda_1)\, . \notag
\end{align}
We also need to define the following spaces
\begin{align}
Q_i:= \prod_{l=1}^{i-1} \mbp \times( \mbp \times V )\times \prod_{l=i}^{k-1} \hat{V}\, , \notag \\
O_i:= \prod_{l=1}^{i-1} \mbp \times (V \times \hat{V} )\times \prod_{l=i}^{k-1} \hat{V}\, , \notag \\
P_i:= \prod_{l=1}^{i-1} \mbp \times (\mbp \times V \times \hat{V}) \times \prod_{l=i}^{k-1} \hat{V}\, , \notag
\end{align}
together with the corresponding projections
\begin{align}
&\pi^{Q_i}: Q_i \lra \mbp \times V\, , \notag \\
&\pi^{O_i}: O_i \lra V \times \hat{V}\, , \notag \\
&\pi^{P_i}: P_i \lra \mbp  \times V \times \hat{V}\, . \notag
\end{align}
This gives rise to the following commutative diagram
$$
\xymatrix{&&& P_i \ar@/_1.3pc/[ddll]_<<<<<<<<<{\pi^{P_i}} \ar[dr]^{q_{23}^i} \ar[d]^{q_{13}^i} \ar[dl]^{q_{12}^i}& \\ & & O_i \ar@/_2.5pc/[ddll]_{\pi^{O_i}}|<<<<<\hole & T_i \ar[ddll]^{\pi^{T_i}}|>>>>>>>>>>>\hole & Q_i \ar@/^2pc/[ddll]^{\pi^{Q_i}} \\ & \mbp \times V \times \hat{V} \ar@/^/[dr]_>>>>>{q_{23}} \ar[d]^{q_{13}} \ar@/_/[dl]^{q_{12}} & & &\\ V \times \hat{V} & \mbp \times \hat{V} & \mbp \times V &}
$$

We consider now the following lifts of the kernels $\widetilde{i}_+\widetilde{S}_u$ on $\mbp \times V$ to $Q_i$:
\[
\widetilde{S}^{i}_1 := (\pi^{Q_i})^+ \widetilde{i}_+\mco_{\widetilde{V}}, \quad \widetilde{S}^{i}_{U^\dag} := (\pi^{Q_i})^+ \widetilde{i}_+ (j_{\dot{V} \dag} \mco_{\dot{V}}),  \quad \widetilde{S}^{i}_\delta := (\pi^{Q_i})^+ \widetilde{i}_+ (i_{E +} \mco_E) %\quad \widetilde{S}^{i}_{1/t} := (\pi^{Q_i})^+ \widetilde{i}_+ (j_{\dot{V} +} \mco_{\dot{V}}),
\]
and
\[
L_i := \pi^{O_i +} L 
\]
on $O_i$. Notice that we have $\widetilde{R}^i_{\hat{u}} \simeq \pi^{T_i +} R'_{\hat{u}}$.
\begin{prop}\label{prop:RadonparttoFLpart}We have the following isomorphisms of kernels in $D^b(\mbp \times \hat{V})$:
\[
\widetilde{S}^{i}_u \diamond L_i \simeq \widetilde{R}^i_{\hat{u}}\quad \text{for} \quad u = 1, \delta, U^\dag  \quad \hat{u}= \delta , 1,  U^+ \, .
\]
\end{prop}
\begin{proof}
\begin{align}
\widetilde{S}^{i}_u \diamond L_i &\simeq q^{i}_{13+}((q_{23}^{i})^+ \widetilde{S}^i_u\overset{L}{\otimes}(q_{12}^{i})^+ L_i ) \notag \\
&\simeq q^{i}_{13+}((q_{23}^{i})^+(\pi^{Q_i})^+  \widetilde{S}_u\overset{L}{\otimes}(q_{12}^{i})^+ (\pi^{O_i})^+ L ) \notag \\
&\simeq q^{i}_{13+}((\pi^{P_i})^+ q_{23}^{+}\, \widetilde{S}_u \overset{L}{\otimes} (\pi^{P_i})^+ q_{12}^{+} L ) \notag \\
&\simeq q^{i}_{13+}(\pi^{P_i})^+ ( q_{23}^{+}\, \widetilde{S}_u \overset{L}{\otimes}  q_{12}^{+} L ) \notag \\
&\simeq (\pi^{T_i})^+ q_{13+}( q_{23}^{+}\, \widetilde{S}_u \overset{L}{\otimes}  q_{12}^{+} L ) \notag \\
&\simeq (\pi^{T_i})^+ (\widetilde{S}_u \diamond L) \notag \\
&\simeq (\pi^{T_i})^+(R_{\hat{u}}) \label{eq:isoAE} \\
&\simeq R^i_{\hat{u}}\, ,  \notag
\end{align}
where \eqref{eq:isoAE} follows from Proposition \ref{prop:isoAE}.
\end{proof}
\noindent Now define the following spaces
\[
N_{i,j,l}:= \mbp^{\times i} \times  V^{\times j} \times \hat{V}^{\times l}\, .
\]
with $i,j,l \in \mbn_0$ and $i+j+l=k$.
Notice that $T_{k-i,i} = N_{k-i,0,i}$.
We define functors
\[
ext^{ijl}_u : D^b_{qc}(\mcd_{N_{i,j,l}}) \ra D^b_{qc}(\mcd_{N_{i-1, j+1,l}}) \quad \text{for}\; u = \delta,1 %,1/t
\]
and
\[
ext^{ijl}_{U^\dag} : D^b_{coh}(\mcd_{N_{i,j,l}}) \ra D^b_{coh}(\mcd_{N_{i-1, j+1,l}})\, ,
\]
which are lifts of the functors
\begin{align}
ext_u: D^b_{qc}(\mcd_\mbp) &\lra D^b_{qc}(\mcd_V)\qquad \qquad \text{for}\; u = \delta,1\, , \notag \\
M &\mapsto ext_u(M) \simeq \pi_{2 +}(\pi_1^+ M \overset{L}{\otimes} \widetilde{S}_u) \notag
\end{align}
resp. $ext_{U^\dag}: D^b_{coh}(\mcd_\mbp) \ra D^b_{coh}(\mcd_V)$. More precisely, let
\begin{align}
\pi_1^{ijl}: \mbp^{\times (i-1)} \times (\mbp \times V) \times  V^{\times j} \times \hat{V}^{\times l} &\ra N_{i,j,l} =\mbp^{\times (i-1)} \times \mbp \times  V^{\times j} \times \hat{V}^{\times l}\, , \label{eq:pieq1} \\
\pi_2^{ijl}: \mbp^{\times (i-1)} \times (\mbp \times V) \times  V^{\times j} \times \hat{V}^{\times l} &\ra N_{i-1,j+1,l} =\mbp^{\times (i-1)} \times V \times  V^{\times j} \times \hat{V}^{\times l}\, , \label{eq:pieq2} \\
\pi^{ijl}: \mbp^{\times (i-1)} \times (\mbp \times V) \times  V^{\times j} \times \hat{V}^{\times l} &\ra \mbp \times V \label{eq:pieq3}
\end{align}
be the canonical projections, then
\[
ext^{ijl}_u(M) := (\pi_{2}^{ijl})_+\left((\pi_1^{ijl})^+ M \overset{L}{\otimes} (\pi^{ijl})^+\widetilde{S}_u\right)
\]
for $M \in D^b_{qc}(\mcd_{N_{i,j,l}})$ and $u= \delta, 1$ resp. $M \in D^b_{coh}(\mcd_{N_{i,j,l}})$ and $u = U^\dag$.\\

We also want to define various partial Fourier-Laplace transforms 
\[
FL_{ijl}: D^b_{qc}(\mcd_{N_{i,j,l}})  \lra D^b_{qc}(\mcd_{N_{i,j-1,l+1}})\, ,
\]
which lift the corresponding functor
\begin{align}
FL: D^b_{qc}(\mcd_V) &\lra  D^b_{qc}(\mcd_{\hat{V}})\, , \notag \\
M &\mapsto FL(M)\, , \notag
\end{align}
i.e. let
\begin{align}
\chi_1^{ijl}:\mbp^{\times i} \times  V^{\times (j-1)} \times (V \times \hat{V}) \times \hat{V}^{\times l} &\lra  N_{i,j,l} =\mbp^{\times i} \times  V^{\times (j-1)} \times V \times \hat{V}^{\times l}\, , \notag \\
\chi_2^{ijl}:\mbp^{\times i} \times  V^{\times (j-1)} \times (V \times \hat{V}) \times \hat{V}^{\times l} &\lra  N_{i,j-1,l+1} =\mbp^{\times i} \times  V^{\times (j-1)} \times \hat{V} \times \hat{V}^{\times l}\, , \notag \\
\chi^{ijl}:\mbp^{\times i} \times  V^{\times (j-1)} \times (V \times \hat{V}) \times \hat{V}^{\times l} &\lra V \times \hat{V} \notag
\end{align}
be the canonical projections, then
\[
FL_{ijl}(M):= (\chi_{2}^{ijl})_+\left((\chi_1^{ijl})^+ M \overset{L}{\otimes} (\chi^{ijl})^+L\right)\, .
\]
Now consider the following diagram of functors (for readability we just wrote the spaces and not the categories on which the functors are defined)
$$
\xymatrix@C=50pt@R=35pt{& & & & & N_{0,0,k}\\
& & & & N_{1,0,k-1} \ar[r]_{ext_u^{1,0,k-1}} \ar[ur]^{\mcr^1_{\hat{u}}} & N_{0,1,k-1} \ar[u]_{FL_{0,1,k-1}} \\
& & & & &\\
& & N_{k-2,0,2}\ar@{.}[uurr] & & &\\
& N_{k-1,0,1} \ar[r]_{ext_u^{k-1,0,1}} \ar[ur]^{\mcr^{k-1}_{\hat{u}}} & N_{k-2,1,1} \ar[u]_{FL_{k-2,1,1}} & & & \\
N_{k,0,0} \ar[r]_{ext_u^{k,0,0}} \ar[ur]^{\mcr^k_{\hat{u}}} & N_{k-1,1,0} \ar[u]_{FL_{k-1,1,0}} \ar[r]_{ext_u^{k-1,1,0}}  & N_{k-2,2,0} \ar@{.}[rrr] \ar[u]_{FL_{k-2,2,0}} & &  & N_{0,k,0}\ar@{.}[uuuu]}
$$

First notice that the upper triangles commute, i.e. for $M \in D^b_{qc}(\mcd_{\mbp^{\times k}})$ and $u= \delta,1$ resp.\\ $M \in D^b_{coh}(\mcd_{\mbp^{\times k}})$ and $u= U^\dag$ we have the following isomorphisms
\[
\widetilde{\mcr}^i_{\hat{u}}(M) \simeq M \diamond \widetilde{R}^i_{\hat{u}} \simeq M \diamond (\widetilde{S}^i_u \diamond L_i) \simeq (M \diamond \widetilde{S}^i_u) \diamond L_i \simeq FL_{i-1,1,k-i}(M \diamond \widetilde{S}^i_u) \simeq FL_{i-1,1,k-i} \circ ext^{i,0,k-i}_u(M)\,.
\]

\begin{lem}\label{lem:FlparttoFLprod} Let $M \in D^b_{qc}(\mcd_{\mbp^{\times k}})$. We have the following isomorphisms
\[
FL_{0,1,k-1} \circ ext_u^{1,0,k-1} \ldots FL_{k-1,1,0} \circ ext_u^{k,0,0}(M) \simeq FL \circ ext_u^{1,k-1,0} \circ \ldots \circ ext_u^{k,0,0}(M)
\]
\end{lem}
\begin{proof}
To prove the proposition we have to show that the following squares commute
\[
\xymatrix@C=50pt@R=35pt{N_{i,j-1,l+1} \ar[r]^{ext^{i,j-1,l+1}_u}& N_{i-1,j,l+1}\\ N_{i,j,l}\ar[u]_{FL_{i,j,l}}\ar[r]_{ext^{i,j,l}_u} & N_{i-1,j+1,l}\ar[u]_{FL_{i-1,j+1,l}}}
\]
for $i,j,l \in \mbn_0$ with $i+j+l=k$.
It is enough to prove the commutativity of the following diagram:
\[
\xymatrix@C=50pt@R=35pt{X \times \mbp \times \hat{V} \ar[r]^{ext_u^2}& X \times V \times \hat{V}\\
X \times \mbp \times V \ar[u]_{FL_2}\ar[r]_{ext_u^1} &  X \times V \times V\ar[u]_{FL_1}}
\]
where $X$ is some smooth algebraic variety and the functors $FL_1, FL_2, ext^1_u, ext^2_u$ are defined below. Consider the following maps
\begin{align}
p_1 : X \times \mbp \times V \times V &\lra X \times \mbp \times V\, , \notag \\
(x,[v^1],v^2,v^3) &\mapsto (x,[v^1],v^2)\, , \notag \\
p_2 : X \times \mbp \times V \times V &\lra X \times V \times V\, , \notag \\
(x,[v^1],v^2,v^3) &\mapsto (x,v^2,v^3)\, , \notag \\
p : X \times \mbp \times V \times V &\lra \mbp \times V\, , \notag \\
(x,[v^1],v^2,v^3) &\mapsto ([v^1],v^3)\, . \notag
\end{align}
The functor $ext_u^1$ is defined by:
\[
ext_u^1(M) := p_{2 +}(p_1^+(M) \overset{L}{\otimes} p^+ (\widetilde{S}_u))\, .
\]
Consider the following maps
\begin{align}
q_1 : X \times \mbp \times V \times \hat{V} &\lra X \times \mbp \times \hat{V}\, , \notag \\
(x,[v^1],v^2,\lambda^1) &\mapsto (x,[v^1],\lambda^1)\, , \notag \\
q_2 : X \times \mbp \times V \times \hat{V} &\lra X \times V \times \hat{V}\, , \notag \\
(x,[v^1],v^2,\lambda^1) &\mapsto (x,v^2,\lambda^1)\, , \notag \\
q : X \times \mbp \times V \times \hat{V} &\lra \mbp \times V\, , \notag \\
(x,[v^1],v^2, \lambda^1) &\mapsto ([v^1],v^2)\, . \notag
\end{align}
The functor $ext_u^2$ is defined by:
\[
ext_u^2(M) = q_{2 +}(q_1^+(M) \overset{L}{\otimes} q^+ (\widetilde{S}_u))\, .
\]
Now define the maps
\begin{align}
\chi_1 : X \times V \times V \times \hat{V} &\lra X \times V \times V\, , \notag \\
(x,v^1,v^2,\lambda^1) &\mapsto (x,v^1,v^2)\, , \notag \\
\chi_2 : X \times V \times V \times \hat{V} &\lra X \times V \times \hat{V}\, , \notag \\
(x,v^1,v^2,\lambda^1) &\mapsto (x,v^2,\lambda^1)\, , \notag \\
\chi : X \times V \times V \times \hat{V} &\lra V \times \hat{V}\, , \notag \\
(x,v^1,v^2, \lambda^1) &\mapsto (v^1,\lambda^1)\, . \notag
\end{align}
The functor $FL_1$ is defined by:
\[
FL_1(M) := \chi_{2 +}(\chi_1^+(M) \overset{L}{\otimes} \chi^+ L)\, .
\]

Now define
\begin{align}
\kappa_1 : X \times \mbp \times V \times \hat{V} &\lra X \times \mbp \times V\, , \notag \\
(x,[v^1],v^2,\lambda^1) &\mapsto (x,[v^1],v^2)\, , \notag \\
\kappa_2 : X \times \mbp \times V \times \hat{V} &\lra X \times \mbp \times \hat{V}\, , \notag \\
(x,[v^1],v^2,\lambda^1) &\mapsto (x,[v^1],\lambda^1)\, , \notag \\
\kappa : X \times \mbp \times V \times \hat{V} &\lra V \times \hat{V}\, , \notag \\
(x,[v^1],v^2, \lambda^1) &\mapsto (v^2,\lambda^1)\, . \notag
\end{align}
The functor $FL_2$ is defined by:
\[
FL_2(M) := \kappa_{2 +}(\kappa_1^+(M) \overset{L}{\otimes} \kappa^+ L)\, .
\]
Consider the following diagram

$$
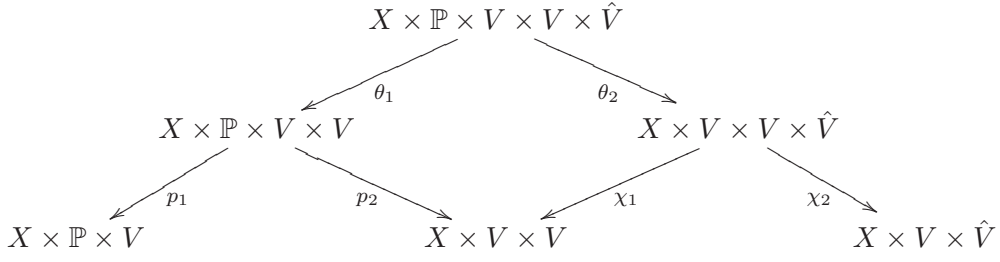
\xymatrix@C=-1pt{ & & X \times \mbp \times V \times V \times \hat{V} \ar[dl]^{\theta_1} \ar[dr]_{\theta_2} & & \\ & X \times \mbp  \times V \times V \ar[dl]^{p_1} \ar[dr]_{p_2} & & X \times V \times V \times \hat{V} \ar[dl]^{\chi_1} \ar[dr]_{\chi_2} & \\ X \times \mbp \times V & & X \times V \times V & & X \times V \times \hat{V}}
$$

where the maps $\theta_1$ and $\theta_2$ are defined as follows
\begin{align}
\theta_1: X \times \mbp \times V \times V \times \hat{V} &\lra X \times \mbp \times V \times V\, , \notag \\
(x,[v^1],v^2,v^3,\lambda^1) &\mapsto (x,[v^1],v^2,v^3)\, , \notag \\
\theta_2: X \times \mbp \times V \times V \times \hat{V} &\lra X \times V \times V \times \hat{V}\, , \notag \\
(x,[v^1],v^2,v^3,\lambda^1) &\mapsto (x,v^2,v^3,\lambda^1)\, . \notag
\end{align}
Notice that the square is a cartesian diagram. We have
\begin{align}
FL_1 \circ ext_u^1(M) &= \chi_{2+}\left( \chi_1^+ p_{2+} \left( p_1^+(M) \overset{L} \otimes p^+ \widetilde{S}_u\right) \overset{L}{\otimes} \chi^+ L\right) \notag \\
&= \chi_{2+} \left( \theta_{2 +} \theta_1^+ \left( p_1^+(M) \overset{L}{\otimes} p^+ \widetilde{S}_u\right) \overset{L}{\otimes} \chi^+L \right)\notag \\
&= \chi_{2+} \left( \theta_{2 +}  \left( \theta_1^+ p_1^+(M) \overset{L}{\otimes} \theta_1^+ p^+ \widetilde{S}_u\right) \overset{L}{\otimes} \chi^+L \right)\notag \\
&= \chi_{2+} \theta_{2 +} \left(   \left( \theta_1^+ p_1^+(M) \overset{L}{\otimes} \theta_1^+p^+ \widetilde{S}_u\right) \overset{L}{\otimes}
\theta_2^+ \chi^+ L \right)\notag \\
&= (\chi_{2}\circ \theta_{2})_+ \left(   \left( (p_1 \circ \theta_1)^+(M) \overset{L}{\otimes} (p \circ \theta_1)^+ \widetilde{S}_u\right) \overset{L}{\otimes}
(\chi \circ \theta_2)^+L \right)\, .\notag
\end{align}
Now consider the other diagram

$$
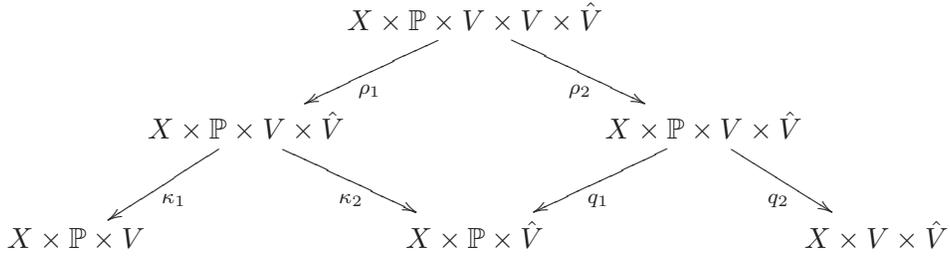
\xymatrix@C=-5pt{ & & X \times \mbp \times V \times V \times \hat{V} \ar[dl]^{\rho_1} \ar[dr]_{\rho_2} & & \\ & X \times \mbp  \times V \times \hat{V} \ar[dl]^{\kappa_1} \ar[dr]_{\kappa_2} & & X \times \mbp \times V \times \hat{V} \ar[dl]^{q_1} \ar[dr]_{q_2} & \\ X \times \mbp \times V & & X \times \mbp \times \hat{V} & & X \times V \times \hat{V}}
$$

where the maps $\rho_1$ and $\rho_2$ are defined as follows:
\begin{align}
\rho_1: X \times \mbp \times V \times V \times \hat{V} &\lra X \times \mbp \times V \times \hat{V}\, , \notag\\
(x,[v^1],v^2,v^3, \lambda^1) &\mapsto (x,[v^1],v^2, \lambda^1)\, , \notag \\
\rho_2: X \times \mbp \times V \times V \times \hat{V} &\lra X \times \mbp \times V \times \hat{V}\, , \notag\\
(x,[v^1],v^2,v^3, \lambda^1) &\mapsto (x,[v^1],v^3, \lambda^1)\, . \notag
\end{align}
This makes the square a cartesian diagram. We have
\begin{align}
ext_u^2 \circ FL_2(M) &= q_{2 +}\left(q_1^+ \kappa_{2+}\left(\kappa_1^+(M) \overset{L}{\otimes} \kappa^+{L}\right) \overset{L}{\otimes} q^+ \widetilde{S}_u \right) \notag \\
&\simeq q_{2 +}\left(\rho_{2+}\rho_1^+\left(\kappa_1^+(M) \overset{L}{\otimes} \kappa^+{L}\right) \overset{L}{\otimes} q^+ \widetilde{S}_u \right) \notag \\
&\simeq q_{2 +}\left(\rho_{2+}\left(\rho_1^+\kappa_1^+(M) \overset{L}{\otimes} \rho_1^+\kappa^+{L}\right) \overset{L}{\otimes} q^+ \widetilde{S}_u \right) \notag \\
&\simeq q_{2 +}\rho_{2+}\left(\left(\rho_1^+\kappa_1^+(M) \overset{L}{\otimes} \rho_1^+\kappa^+{L}\right) \overset{L}{\otimes} \rho_{2}^+q^+ \widetilde{S}_u \right) \notag \\
&\simeq (q_{2}\circ \rho_{2})_+ \left(\left((\kappa_1\circ \rho_1)^+(M) \overset{L}{\otimes} (\kappa \circ \rho_1)^+{L}\right) \overset{L}{\otimes} (q \circ \rho_{2})^+ \widetilde{S}_u \right)\, . \notag
\end{align}
Notice that we have $q_2 \circ \rho_2 = \chi_2 \circ \theta_2$, $\kappa_1 \circ \rho_1 = p_1 \circ \theta_1$, $\kappa \circ \rho_1 = \chi \circ \theta_2$ and $q \circ \rho_2 =  p \circ \theta_1$. This, together with the associativity of $\overset{L}{\otimes}$ shows
\[
ext^2_u \circ FL_2 = FL_1 \circ ext^1_u\, .
\]
But as observed above, this shows the claim.
\end{proof}

Consider the following diagram, which is the $k$-th product of diagram \eqref{eq:diagram} above:
\[
\xymatrix{ && \dot{V}^{\times k} \ar[d]^{\prod j_{\dot{V}}}\ar[drr]^{\prod j} \ar[dll]_{\prod \pi} && \\
\mbp^{\times k} && \widetilde{V}^{\times k} \ar[rr]^{\prod \widetilde{j}} \ar[ll]_{\prod \widetilde{\pi}} && V^{\times k}\, . \\
&&  E^{\times k} \ar[u]_{\prod i_{E}} \ar[urr]_{\prod j_{0}} \ar[ull]^{\prod \pi^E} &&}
\]

Define the following functors from $D^b_{qc}(\mcd_{\mbp^{\times k}})$ to $D^b_{qc}(\mcd_{V^{\times k}})$:
\begin{align}
ext^{\prod}_\delta(M) &:= (\Uppi\, j_{0})_+ (\Uppi\, \pi^E)^+ M  \notag \\
ext_1^{\prod}(M) &:= (\Uppi\, \widetilde{j})_+ (\Uppi\, \widetilde{\pi})^+ M \notag \\
%ext^{\prod}_{1/t}(M) &:= (\Uppi\, j)_+ (\Uppi\, \pi)^+ M  \notag
\end{align}
resp. 
\[
ext^{\prod}_{U^\dag}(M) := (\Uppi\, j)_\dag (\Uppi\, \pi)^+ M 
\]
from $D^b_{coh}(\mcd_{\mbp^{\times k}})$ to $D^b_{coh}(\mcd_{V^{\times k}})$.

\begin{lem}\label{lem:extprodtodiagext}
Let $M \in D^b_{qc}(\mcd_{\mbp^{\times k}})$ and $u=\delta,1$ resp. $M \in D^b_{coh}(\mcd_{\mbp^{\times k}})$ and $u =U^\dag$. We have the following isomorphisms
\[
ext^{1,k-1,0}_u \circ \ldots \circ ext^{k,0,0}_u  (M) \simeq ext^{\prod}_u (M)\, . 
\]
where $ext^{\prod}_u$ are the extension functors defined above.
\end{lem}
\begin{proof}
We begin by proving the case $u=1$. Set $M_{i,j} := \mbp^{\times i} \times \widetilde{V}^{\times (k-i-j)} \times V^{\times j}$ for $i+j =k$ and define morphisms
\begin{align}
\sigma_1^{ij} := id_{\mbp}^{\times i} \times \widetilde{\pi} \times id_{\widetilde{V}}^{\times k-i-j-1} \times id_V^{\times j}: M_{i,j} &\lra M_{i+1,j}\, , \notag \\
\sigma_2^{ij} := id_{\mbp}^{\times i} \times id_{\widetilde{V}}^{\times k-i-j-1}\times \widetilde{j} \times id_V^{\times j}: M_{i,j} &\lra M_{i,j+1}\, , \notag
\end{align}
where $\widetilde{\pi}: \widetilde{V} \ra \mbp$ and $\widetilde{j} :\widetilde{V} \ra V$ are the maps from equation \eqref{eq:blowupmaps}.\\

Consider the following diagram
$$
\xymatrix{M_{0,0} \ar[r]_{\sigma^{0,0}_2} \ar[d]_{\sigma^{0,0}_1} & M_{0,1} \ar[d]_{\sigma^{0,1}_1} \ar@{.}[rr] & & M_{0,k-1}\ar[r]_{\sigma^{0,k-1}_2} \ar[d]_{\sigma^{0,k-1}_1} & M_{0,k} \\
M_{1,0} \ar[r]_{\sigma^{1,0}_2} \ar@{.}[dd]  & M_{1,1}\ar@{.}[dd]\ar@{.}[rr]  & & M_{1,k-1} &\\
&&&&&\\
M_{k-1,0} \ar[d]_{\sigma^{k-1,0}_1} \ar[r]_{\sigma^{k-1,0}_2}& M_{k-1,1}\ar@{.}[uurr]&&&&\\
M_{k,0} &&&&&}
$$
where the squares are cartesian diagrams. Notice that $M_{k-i,i}= N_{k-i,i,0} = \mbp^{\times (k-i)} \times V^{\times i}$ and 
\[
ext^{k-i,i,0}_1(M) \simeq (\sigma_{2}^{k-i-1,i})_+ (\sigma_1^{k-i-1,i})^+(M)
\]
for $i \in \{0,\ldots,k-1\}$ (cf. \cite[lemma 1]{AE}). We therefore conclude by base change that
\begin{align}
ext^{1,k-1,0}_1 \ldots \circ ext^{k,0,0}_1(M) &\simeq (\sigma_{2}^{0,k-1})_+\ldots (\sigma_{2}^{0,0})_+(\sigma_1^{0,0})^+ \ldots (\sigma_1^{k-1,0})^+(M) \notag \\
&\simeq  (\sigma_{2}^{0,k-1} \circ \ldots \circ \sigma_{2}^{0,0})_+(\sigma_1^{k-1,0} \circ \ldots \circ \sigma_1^{0,0})^+(M)\, . \notag 
\end{align}
Notice that $\sigma_1^{k-1,0} \circ \ldots \circ \sigma_1^{0,0} = \Uppi\,  \widetilde{\pi} :\widetilde{V}^{\times k}  \ra \mbp^{\times k}$ and $\sigma_{2}^{0,k-1} \circ \ldots \circ \sigma_{2}^{0,0} = \Uppi\, \widetilde{j} : \widetilde{V}^{\times k} \ra V^{\times k}$.
This shows the claim for $u=1$.

The proof for $u= U^\dag$ is virtually the same if one defines $M_{i,j}:= \mbp^{ \times i} \times \dot{V}^{\times (k-i-j)} \times V^{\times j}$ with
\begin{align}
\sigma_1^{ij} := id_{\mbp}^{\times i} \times \pi \times id_{\dot{V}}^{\times k-i-j-1} \times id_V^{\times j}: M_{i,j} &\lra M_{i+1,j}\, , \notag \\
\sigma_2^{ij} := id_{\mbp}^{\times i} \times id_{\dot{V}}^{\times k-i-j-1}\times j \times id_V^{\times j}: M_{i,j} &\lra M_{i,j+1}\, , \notag
\end{align}
where $\pi:\dot{V} \ra \mbp$ is the canonical projection and $j: \dot{V} \ra V$ the natural inclusion.
In this case
\begin{align}
%ext^{k-i,i,0}_{1/t}(M) &\simeq (\sigma_{2}^{k-i-1,i})_+ (\sigma_1^{k-i-1,i})^+(M)\, , \notag \\
ext^{k-i,i,0}_{U^\dag}(M) &\simeq (\sigma_{2}^{k-i-1,i})_\dag (\sigma_1^{k-i-1,i})^+(M)\notag
\end{align}
as well as $\sigma_1^{k-1,0} \circ \ldots \circ \sigma_1^{0,0} = \Uppi\, \pi :\dot{V}^{\times k}  \ra \mbp^{\times k}$ and $\sigma_{2}^{0,k-1} \circ \ldots \circ \sigma_{2}^{0,0} = \Uppi\, j : \dot{V}^{\times k} \ra V^{\times k}$.

It remains to prove the case $u = \delta$. Define $L_{i,j} := \mbp^{\times i} \times \{0\}^{k-i-j} \times V^{\times j}$ for $i+j \leq k$ and define morphisms
\begin{align}
\alpha_1^{ij} := id^{\times i-1}_{\mbp} \times p \times id_{\{ 0\}}^{\times k-i-j} \times id^{\times j}_V : L_{i,j} &\lra L_{i-1,j}\, , \notag \\
\alpha_2^{ij} := id^{\times i}_{\mbp} \times id_{\{0\}}^{\times k-i-j-1} \times i_0 \times id^{\times j}_V: L_{i,j} &\lra L_{i,j+1}\, , \notag
\end{align}
where $p: \mbp \ra \{0\}$ is the map to the point and $i_0 : \{0\} \ra V$ is the canonical inclusion.

Consider the following diagram
$$
\xymatrix{L_{0,0}\ar[r]^{\alpha_2^{0,0}}  & L_{0,1}   \ar@{.}[rr] & & L_{0,k-1} \ar[r]^{\alpha_2^{0,k-1}}& L_{0,k}  \\
L_{1,0} \ar[u]^{\alpha_1^{1,0}}  \ar@{.}[dd] \ar[r]^{\alpha_2^{1,0}}  & L_{1,1}\ar@{.}[dd]\ar@{.}[rr]  \ar[u]^{\alpha_1^{1,1}}  & & L_{1,k-1} \ar[u]^{\alpha_1^{1,k-1}} &\\
&&&&&\\
L_{k-1,0} \ar[r]^{\alpha_2^{k-1,0}}& L_{k-1,1}\ar@{.}[uurr]&&&&\\
L_{k,0} \ar[u]^{\alpha_1^{k,0}}&&&&&}
$$
where the squares are commutative diagrams. For the following isomorphisms, notice that $L_{k-i,i}=N_{k-i,i,0}= \mbp^{\times (k-i)} \times V^{\times i}$, and recall the following maps
\begin{align}
\pi_1^{k-i,i,0}: N_{k-i,i+1,0} &\lra N_{k-i,i,0} = L_{k-i,i}\, , \notag \\
\pi_2^{k-i,i,0}: N_{k-i,i+1,0} &\lra N_{k-i-1,i+1,0} = L_{k-i-1,i+1}\, , \notag \\
\pi^{k-i,i,0}: N_{k-i,i+1,0} &\lra \mbp \times V\, , \notag
\end{align}
which were defined in equation \eqref{eq:pieq1}, \eqref{eq:pieq2} and \eqref{eq:pieq3}.

 and denote by $\xi_i : L_{k-i,i} \ra  N_{k-i,i+1,0}$ the embedding with image $\mbp^{\times (k-i)}\times \{0\} \times V^{\times i}$. 
\begin{align}
ext^{k-i,i,0}_\delta(M) &\simeq (\pi_2^{k-i,i,0})_+\left( (\pi_1^{k-i,i,0})^+ M \overset{L}{\otimes} (\pi^{k-i,i,0})^+ \widetilde{i}_+\widetilde{S}_\delta \right) \notag \\
&\simeq (\pi_2^{k-i,i,0})_+\left( (\pi_1^{k-i,i,0})^+ M \overset{L}{\otimes} \xi_{i+} \mco_{L_{k-i,i}} \right) \notag \\
&\simeq (\pi_2^{k-i,i,0})_+ \xi_{i+} \left( \xi_i^+ (\pi_1^{k-i,i,0})^+ M \overset{L}{\otimes}  \mco_{L_{k-i,i}} \right) \notag \\
&\simeq (\pi_2^{k-i,i,0})_+ \xi_{i+} M \notag \\
&\simeq (\alpha_2^{k-i-1,i})_+ \alpha_{1+}^{k-i,i} (M)\,.
\end{align}
for $i \in \{0,\ldots,k-1\}$, where the second isomorphism follows from the fact that the image of the exceptional divisor $E \subset \widetilde{V}$ under the embedding $\widetilde{i}: \widetilde{V} \ra \mbp \times V$ is equal to $\mbp \times \{0\}$ and the fourth isomorphism follows from the fact that $\xi_i$ is a section of $\pi_1^{k-i,i,0}$.  We therefore conclude by the commutativity of the diagram that
\begin{align}
ext_\delta^{1,k-1,0} \circ \ldots \circ ext_\delta^{k,0,0}(M) &\simeq (\alpha_2^{0,k-1})_+ \ldots (\alpha_2^{0,0})_+ \alpha_{1+}^{1,0} \ldots \alpha_{1+}^{k,0}(M) \notag \\
&\simeq (\alpha_2^{0,k-1} \circ \ldots \circ \alpha_2^{0,0})_+ (\alpha_1^{1,0} \circ \ldots \circ \alpha_1^{k,0})_+ M \notag  \\
&\simeq (i_{0})_+ p_+ M \notag \\
&\simeq (\Uppi\, j_{0})_+ (\Uppi\, \pi^E)^+ M \notag \\
&= ext^{\prod}_\delta(M)\, , \notag
\end{align}
where $p: \mbp^{\times k} \ra \{0\}$ is the map to a point and $i_{0}: \{0\} \ra V^{\times k}$ is the embedding with image $\{0\} \times \ldots \times \{0\}$ and $\Uppi \, j_{0}$ resp. $\Uppi \, \pi^E$ were defined above.
\end{proof}

Finally, we have to compare the functors $ext^{\prod}_u$ and $ext^k_u$.
\begin{lem}\label{lem:compext}
Let $M \in D^b_{qc}(\mcd_{\mbp})$ and $u=\delta,1$ resp. $M \in D^b_{coh}(\mcd_{\mbp})$ and $u =U^\dag$. We have the following isomorphisms
\[
ext^{\prod}_u \circ \Delta_+  (M) \simeq ext^{k}_u (M)\, . 
\]
where $ext^{k}_u$ are the extension functors defined in Section \ref{sec:statresult}.
\end{lem}
\begin{proof}
In order to prove the lemma, consider the following diagrams
\[
\xymatrix{ & V^{\times k} & \\ \mbv \ar[ur]^{j_k} \ar[rr]^{\widetilde{i}_k} \ar[d]^{\pi_k}& & \widetilde{V}^{\times k} \ar[ul]_{\prod \widetilde{j}} \ar[d]^{\prod \widetilde{\pi}}\\ \mbp \ar[rr]^{\Delta} && \mbp^{\times k}}\quad
\xymatrix{ & V^{\times k} & \\ \dot{\mbv} \ar[ur]^{j^\ast} \ar[rr] \ar[d]^{\pi^\ast}& & \dot{V}^{\times k} \ar[ul]_{\prod j} \ar[d]^{\prod \pi}\\ \mbp \ar[rr]^{\Delta} && \mbp^{\times k}} \quad
\xymatrix{ & V^{\times k} & \\ E \ar[ur]^{j_E} \ar[rr] \ar[d]^{\pi^E}& & E^{\times k} \ar[ul]_{\prod j_0} \ar[d]^{\prod \pi^E}\\ \mbp \ar[rr]^{\Delta} && \mbp^{\times k}}  
\]
where the lower squares are cartesian and the upper triangles are commutative, respectively. We will prove the lemma in the case $u=1$, the other cases are similar:
\[
ext^{\prod}_1 \circ \Delta_+(M)\; = \; (\Uppi\, \widetilde{j})_+ (\Uppi\, \widetilde{\pi})^+ \Delta_+ (M)\;  \simeq\; (\Uppi\, \widetilde{j})_+( \widetilde{i}_k)_+ \widetilde{\pi}_k^+(M)\; \simeq\;  (\widetilde{j}_k)_+ \widetilde{\pi}_k^+ (M)\; = \; ext^k_1(M)
\]
\end{proof}

We are now able to give the proof of Proposition \ref{prop:higheraffRad}.

\begin{proof}[Proof of proposition \ref{prop:higheraffRad}]
For $u = \delta,1 , U^\dag$ and $\hat{u} = 1, \delta , U^+$ we have the following isomorphisms:
\begin{align}
\mcr'_{\hat{u}}(M) &\simeq \widetilde{\mcr}_{\hat{u}}(\Delta_+(M)) \notag \\
&\simeq \widetilde{\mcr}^1_{\hat{u}} \circ \ldots \circ \widetilde{\mcr}^k_{\hat{u}} (\Delta_+ (M)) \notag \\
&\simeq FL_{0,1,k-1} \circ ext_u^{1,0,k-1} \ldots FL_{k-1,1,0} \circ ext_u^{k,0,0}(\Delta_+(M)) \notag \\
&\simeq FL \circ ext_u^{1,k-1,0} \circ \ldots \circ ext_u^{k,0,0}(\Delta_+(M)) \notag \\
&\simeq FL \circ ext^{\prod}_u (\Delta_+ (M))\, , \notag \\
&\simeq FL \circ ext^k_u (M)\, , \notag
\end{align}

where the first isomorphism follows from Lemma \ref{lem:RadontoRadonprod}, the second isomorphism from Lemma \ref{lem:RadonprodtoRadonpart}, the third isomorphism follows form Proposition \ref{prop:RadonparttoFLpart}, the fourth isomorphism follows from Lemma \ref{lem:FlparttoFLprod}, the fifth isomorphism follows from Lemma \ref{lem:extprodtodiagext} and the last isomorphism from Lemma \ref{lem:compext}.

This gives rise to the following triangles
\[
\xymatrix{\mcr'_1(M) \ar[r]  \ar[d]^\simeq & \mcr'_{1/t}(M) \ar[r] \ar@{.>}[d]  & \mcr'_{\delta}(M) \ar[r]^{+1} \ar[d]^\simeq  & \\
FL \circ ext^k_\delta(M) \ar[r] & FL \circ ext^k_{Y}(M) \ar[r] & FL \circ ext^k_1(M) \ar[r]^<<<<{+1} & }
\]
\[
\xymatrix{\mcr'_\delta(M) \ar[r]  \ar[d]^\simeq & \mcr'_{Y}(M) \ar[r] \ar@{.>}[d]  & \mcr'_{1}(M) \ar[r]^{+1} \ar[d]^\simeq  & \\
FL \circ ext_1(M) \ar[r] & FL \circ ext_{1/t}(M) \ar[r] & FL \circ ext_\delta(M) \ar[r]^<<<<{+1} & }
\]
\[
\xymatrix{\mcr'_{A^+}(M) \ar[r]  \ar@{.>}[d] & \mcr'_{1}(M) \ar[r] \ar[d]^\simeq  & \mcr'_{U^+}(M) \ar[r]^{+1} \ar[d]^\simeq  & \\
FL \circ ext_{A^\dag}(M) \ar[r] & FL \circ ext_{\delta}(M) \ar[r] & FL \circ ext_{U^\dag}(M) \ar[r]^<<<<{+1} & }
\]
Notice that that the connecting homomorphisms in the first two diagrams resp. the second square in the last diagram commutes by the functoriality of the proof. Hence the third arrow exists in each diagram and is an isomorphism in $D^b_{coh}(\mcd_{\hat{W}})$.
\end{proof}

We can now give the proof of Theorem \ref{thm:higherRad}.

\begin{proof}[Proof of Theorem \ref{thm:higherRad}]
Using Proposition \ref{prop:higheraffRad} it remains to prove
\[
r^+\mcr'_{\hat{u}}(M) \simeq \pi^+ \mcr_{\hat{u}}(M)
\]
Let $Z'':= \{[v], \lambda^1, \ldots, \lambda^k \in \mbp \times S(k,n) \mid \lambda^1(v) = \ldots = \lambda^k(v) = 0 \}$  and $C'' := \{[v], \lambda^1, \ldots, \lambda^k \in \mbp \times S(k,n) \mid \exists l \in \{1, \ldots ,k \} \; \text{with}\; \lambda^l(v) \neq 0 \}$.

\noindent We will prove the case $r^+\mcr'_{\delta}(M) \simeq \pi^+ \mcr_{\delta}(M)$ the other cases are similar or simpler.

\noindent Consider the following commutative diagram
$$
\xymatrix{ && Z' \ar[dll]_{\pi_1^{Z'}} \ar[rr]^{\pi_2^{Z'}} && \hat{W}\\ \mbp && Z'' \ar[d]^{\pi_{Z''}} \ar[u]_{r_{Z''}} \ar[ll]_{\pi_1^{Z''}} \ar[rr]^{\pi_2^{Z''}} && S(k,n)\ar[u]_r \ar[d]^\pi \\ && Z \ar[rr]_{\pi_2^Z} \ar[ull]^{\pi_1^Z} && \mbg}
$$
where $\pi_1^{Z''}$ resp. $\pi_2^{Z''}$ are the canonical projections restricted to $Z''$ and $r_{Z''}$ resp. $\pi_{Z''}$ were chosen such that the upper resp. lower square becomes cartesian. One has
\begin{align}
r^+\mcr'_\delta(M) &= r^+ \pi_{2+}^{Z'} (\pi_{1}^{Z'})^+ M \notag \\
&\simeq \pi_{2 +}^{Z''} r_{Z''}^+ (\pi_1^{Z'})^+ M \notag \\
&\simeq \pi_{2 +}^{Z''} (\pi_1^{Z''})^+ M \notag \\
&\simeq \pi_{2 +}^{Z''} \pi_{Z''}^+ (\pi_1^Z)^+ M \notag \\
&\simeq \pi^+ \pi_{2+}^Z (\pi_1^Z)^+ M \notag \\
&= \pi^+ \mcr_\delta M \notag
\end{align}
\end{proof}
\bibliographystyle{amsalpha}
\bibliography{radon}
\vspace*{1cm}

\noindent
Thomas Reichelt\\
Mathematisches Institut \\
Universit\"at Heidelberg\\
69120 Heidelberg\\
Germany\\
treichelt@mathi.uni-heidelberg.de
\end{document}